\documentclass[10pt]{article}

\DeclareMathAlphabet{\mathpzc}{OT1}{pzc}{m}{it}

\usepackage{amsmath}
\usepackage{amssymb}
\usepackage{amsfonts}
\usepackage{amsthm}
\usepackage{mathrsfs}
\usepackage{ifsym}
\usepackage{fontenc}
\usepackage{textcomp}
\usepackage{tipa}
\usepackage{enumerate}
\usepackage[pdftex]{color, graphicx}

\begin{document}
 
\title{{\bf  On the Hasse invariants of the Tate normal forms $E_5$ and $E_7$} \footnote{{\it 2020 Mathematics Subject Classification}. 14H52, 11G15, 11G20, 11R29, 11R37} \footnote{{\it Key words and phrases}. Tate normal form, Hasse invariant, class number, class fields, Rogers-Ramanujan continued fraction, Fricke group}}       
\author{Patrick Morton}        
\date{}          
\maketitle

\begin{abstract}  A formula is proved for the number of linear factors over $\mathbb{F}_l$ of the Hasse invariant of the Tate normal form $E_5(b)$ for a point of order $5$, as a polynomial in the parameter $b$, in terms of the class number of the imaginary quadratic field $K=\mathbb{Q}(\sqrt{-l})$, proving a conjecture of the author from 2005.  A similar theorem is proved for quadratic factors with constant term $-1$, and a theorem is stated for the number of quartic factors of a specific form in terms of the class number of $\mathbb{Q}(\sqrt{-5l})$.  These results are shown to imply a recent conjecture of Nakaya on the number of linear factors over $\mathbb{F}_l$ of the supersingular polynomial $ss_l^{(5*)}(X)$ corresponding to the Fricke group $\Gamma_0^*(5)$.  The degrees and forms of the irreducible factors of the Hasse invariant of the Tate normal form $E_7$ for a point of order $7$ are determined, which is used to show that the polynomial $ss_l^{(N*)}(X)$ for the group $\Gamma_0^*(N)$ has roots in $\mathbb{F}_{l^2}$, for any prime $l  \neq N$, when $N \in \{2,3,5,7\}$.
\end{abstract}

\section{Introduction.}

The Tate normal form of an elliptic curve with a point of order $5$ is the curve
$$E_5: Y^2+(1+b)XY+bY=X^3+bX^2.$$
In \cite{mor} I proved several properties of the Hasse invariant of this curve, which for a given prime $l \ge 3$ is the polynomial
\begin{align*}
\hat H_{5,l}(b) & = (b^4+12b^3+14b^2-12b+1)^r(b^2+1)^s(b^4+18b^3+74b^2-18b+1)^s\\
& \times b^{5n_l}(1-11b-b^2)^{n_l}J_l\left(\frac{(b^4+12b^3+14b^2-12b+1)^3}{b^5(1-11b-b^2)}\right),
\end{align*}
where $n_l=[l/12], r = \frac{1}{2}(1-(-3/l)), s=\frac{1}{2}(1-(-4/l))$; and  $J_l(t)$ is the polynomial
$$J_l(t)\equiv \sum_{k=0}^{n_l}{{2n_l + s\atopwithdelims ( ) 2k + s}{2n_l-2k \atopwithdelims ( ) n_l-k}(-432)^{n_l-k}(t-1728)^k} \ \ (\textrm{mod} \ l).$$
The factors of the Hasse invariant $\hat H_{5,l}(b)$ over the field $\mathbb{F}_l$ are linear, quadratic, or quartic in $b$ (see \cite{mor}, Theorem 6.2).  Linear factors occur only for primes $l \equiv 4$ (mod $5$), and in \cite{mor} I made a conjecture about the number of distinct linear factors in this case.  In this paper I will prove this conjecture, which can be stated as follows. \bigskip

\noindent {\bf Theorem 1.1.} {\it If $l \equiv 4$ (mod $5$) is a prime and $h(-l)$ denotes the class number of the field $\mathbb{Q}(\sqrt{-l})$, then the number of distinct linear factors of the Hasse invariant $\hat H_{5,l}(b)$ is}
$$N_l=\begin{cases}
2h(-l), & \textit{if} \ l \equiv 9 \ (\textit{mod} \ 20);\\
2(3-\left(\frac{2}{l}\right))h(-l), & \textit{if} \ l \equiv 19 \ (\textit{mod} \ 20).
\end{cases}$$ \medskip

It is a consequence of a theorem of Deuring that, for a given prime $l$, the number $N_l$ in Theorem 1.1 is exactly $4$ times the number of supersingular invariants which lie in $\mathbb{F}_l$, which number is given by the formula
$$S(\mathbb{F}_l)=\begin{cases}
   \frac{1}{2} h(-l)   & \text{if} \ l \equiv 1 \ \text{mod} \ 4, \\
   2h(-l)   & \text{if} \ l \equiv 3 \ \text{mod} \ 8, \\
   h(-l) & \text{if} \ l \equiv 7 \ \text{mod} \ 8.
\end{cases}$$
(See \cite[p. 97]{brm}.)  To prove the theorem, I will show the following.  Set
\begin{align*}
f_0(x) & = x^4+12x^3+14x^2-12x+1, \\
f_{1728}(x) & = x^4+18x^3+74x^2-18x+1, \\
f_j(x) & = (x^4+12x^3+14x^2-12x+1)^3 - jx^5(1-11x-x^2),
\end{align*}
where $j$ runs over the supersingular invariants different from $0, 1728$ which lie in $\mathbb{F}_l$.  These latter values are exactly the roots of the polynomial $J_l(t)$ in $\mathbb{F}_l$.  Excepting the factor $b^2+1$ (when $l \equiv 3$ mod $4$), the polynomials $f_0(b), f_{1728}(b), f_j(b)$ make up the individual factors of $\hat H_{5,l}(b)$ which correspond to $j$-invariants in $\mathbb{F}_l$.  Hence, Theorem 1.1 follows from:  \bigskip

\noindent {\bf Theorem 1.2.} {\it If $0, 1728$, or $j \neq 0, 1728 \in \mathbb{F}_l$ is a supersingular $j$-invariant for the prime $l \equiv 4$ (mod $5$), then the corresponding polynomial $f_0(x), f_{1728}(x)$, or $f_j(x)$ has exactly $4$ distinct linear factors over $\mathbb{F}_l$.}  \bigskip

The polynomials $f_j(x)$ (including $j=0$ or $j=1728$) have no repeated factors in $\mathbb{F}_l$, as was noted in \cite[p. 268]{mor}.  Also, $f_j(x) \equiv 0$ (mod $l$) can only have a root in $\mathbb{F}_l$ if $j$ also lies in this field. \medskip

The key to proving Theorem 1.2 lies in considering the factorization of the prime $l$ in certain class fields over the imaginary quadratic field $K=\mathbb{Q}(\sqrt{-l})$.  Since $l \equiv 4$ (mod $5$), the prime $5$ splits in the field $K$.  Denoting the conjugate prime divisors of $(5)$ in the ring of integers $R_K$ of $K$ by $\wp_5, \wp_5'$, the relevant class fields are $L=\Sigma_{\mathfrak{f}}$ and $L=\Sigma_{\mathfrak{f}}\Omega_2$ (when $l \equiv 3$ mod $4$), where $\Sigma_{\mathfrak{f}}$ is the ray class field of conductor $\mathfrak{f}=\wp_5'$ or $\wp_5$ over $K$ and $\Omega_2$ is the ring class field of conductor $2$ over $K$.  (The latter field coincides with the Hilbert class field when $l \equiv 7$ (mod $8$).)  \medskip

The field $K$ is useful in the proof of Theorem 1.2 because of the following considerations.  Using Deuring's theorem that $\sqrt{-l}$ injects into the endomorphism ring of any supersingular elliptic curve $E$ with $j(E) \in \mathbb{F}_l$, it follows that every supersingular $j$-invariant lying in $\mathbb{F}_l$ is a root of one of the class equations $H_{-l}(x)$ or $H_{-4l}(x)$ (mod $l$), where the polynomial $H_{-l}(x)$ occurs only when $l \equiv 3$ (mod $4$).  In fact, the following congruences hold.  If $l \equiv 3 \ (\textrm{mod} \ 4)$ then
\begin{equation}
H_{-l}(x) \equiv (x-1728)R(x)^2 \ \ \textrm{and} \ \ H_{-4l}(x) \equiv (x-1728)S(x)^2 \ (\textrm{mod} \ l),\tag{1.1}
\label{eq:1.1}
\end{equation}
where $R(x)$ and $S(x)$ are products of distinct linear factors which are relatively prime to each other.  Moreover,
\begin{align*}
R(x) \equiv & \ \textrm{gcd}(J_l(x),(x-1728)^{(l-1)/2}-1) \ (\textrm{mod} \ l),\\
S(x) \equiv & \ \textrm{gcd}(xJ_l(x),(x-1728)^{(l-1)/2}+1)  \ (\textrm{mod} \ l).
\end{align*}
(In particular, this shows that $x=0$ is never a zero of $H_{-l}(x)$ (mod $l$), since $x=0$ is not a zero of $J_l(x)$.) On the other hand, if $l \equiv 1$ (mod $4$), then 
$$H_{-4l}(x) \equiv T(x)^2, \ \ T(x) \equiv \ \textrm{gcd}(xJ_l(x),(x-1728)^{(l-1)/2}+1)  \ (\textrm{mod} \ l).\eqno{(1.2)}$$
Note also that $x^r(x-1728)^sJ_l(x) =ss_l(x)$ is the supersingular polynomial in characteristic $l$. See Elkies \cite{el}, Kaneko \cite{ka}, and Brillhart and Morton \cite[Propositions 9, 11]{brm}.\medskip

As is shown in \cite{mor2}, the class fields $\Sigma_{\mathfrak{f}}$ and $\Sigma_{\mathfrak{f}}\Omega_2$ are generated over the field $K=\mathbb{Q}(\sqrt{-l})$ by certain values of the Rogers-Ramanujan continued fraction
\begin{align*}
r(\tau) &= \frac{q^{1/5}}{1+\frac{q}{1+\frac{q^2}{1+ \frac{q^3}{1+\cdots}}}}=\frac{q^{1/5}}{1+} \ \frac{q}{1+} \ \frac{q^2}{1+} \ \frac{q^3}{1+} \dots,\\
& = q^{1/5} \prod_{n \ge 1}{(1-q^n)^{(n/5)}},  \ \ q = e^{2 \pi i \tau}, \ \ \tau \in \mathbb{H}.
\end{align*}
(See \cite{anb}, \cite{ber}, \cite{du}.)  Namely, if $d=lf^2$, then $\Sigma_{\mathfrak{\wp_5'}}\Omega_f=\mathbb{Q}(r(w/5))$ and $\Sigma_{\mathfrak{\wp_5}}\Omega_f=\mathbb{Q}(\zeta^j r(-1/w))$, for a unique primitive $5$-th root of unity $\zeta^j$ and an algebraic integer $w=(v+\sqrt{-d})/2$ in $K$ with $\wp_5^2 \mid w$ and $(N(w),f) = 1$.  The numbers $\eta=r(w/5)$ and $\xi = \zeta^j r(-1/w)$ are conjugate units in these fields and roots of a polynomial $p_d(x) \in \mathbb{Z}[x]$ of degree $4\textsf{h}(-d)$, where $\textsf{h}(-d)$ is the class number of the order $\mathcal{O}=\textsf{R}_{-d}$ of discriminant $-d$ in $K$.  Using the above factorizations modulo $l$, the proof of Theorem 1.2 shows that four roots of the polynomial $f_j(x)$ in $\mathbb{F}_l$ (for $j \neq 0, 1728$) are provided by rational expressions of roots of $p_d(x)$ in $\mathbb{F}_l$, for $-d = -l$ or $-4l$; and the latter choice of $-l$ versus $-4l$ depends only on $j$.  The fact that these are the only four roots of $f_j(x)$ in $\mathbb{F}_l$ follows from calculations using the normal form of the icosahedral group over $\mathbb{F}_l(\zeta)$, where $\zeta$ is a $5$-th root of unity.  (See \cite[Vol. II]{fr1}.)  \medskip

As a corollary of the proof of Theorem 1.2, I show in Section 4 that $p_d(x)$ factors (mod $l$) as a product of the squares of distinct linear factors, when $l \equiv 4$ (mod $5$) and $d=l$ or $4l$, echoing the factorizations of the class equations $H_{-d}(x)$ above.  This implies that the powers of $\eta=r(w/5)$ form an $\mathfrak{l}$-integral basis of the field $\Sigma_{\mathfrak{\wp_5'}}\Omega_f$ over $K$, where $\mathfrak{l}=(\sqrt{-l})$ in $R_K$ and $d=lf^2$ ($f=1$ or $2$).  The latter result is also true when $l \equiv 1$ (mod $5$). \medskip

Using a similar approach, the following theorem is proved in Section 4, regarding irreducible quadratic factors of $\hat H_{5,l}(x)$ with constant term $-1$. \bigskip

\noindent {\bf Theorem 1.3.} {\it Let $l > 5$ be a prime.} \begin{enumerate}[a)]
\item {\it If $l \equiv 4$ mod $5$, the polynomial $\hat H_{5,l}(x)$ has no irreducible factors of the form $q(x) = x^2+ax-1$ over $\mathbb{F}_l$.}
\item {\it If $l \equiv 1$ mod $5$, the number of irreducible quadratic factors $q(x)=x^2+ax-1$ dividing $\hat H_{5,l}(x)$ over $\mathbb{F}_l$ is $2S(\mathbb{F}_l)$, i.e. given by:}
\begin{enumerate}[i)]
\item {\it $h(-l)$, if $l \equiv 1$ mod $4$;}
\item {\it $4h(-l)$, if $l \equiv 3$ mod $8$;}
\item {\it $2h(-l)$, if $l \equiv 7$ mod $8$.}
\end{enumerate}
\end{enumerate}

Further, in Section 7, I will use the following theorem regarding quartic factors of $\hat H_{5,l}(x)$, which only occur when $l \equiv 2,3$ (mod $5$), and certain quadratic factors, which only appear when $l \equiv 1, 4$ (mod $5$).  For this theorem, set $\varepsilon = \frac{-1+\sqrt{5}}{2}$ and $\bar \varepsilon = \frac{-1-\sqrt{5}}{2}$, and let $h(-5l)$ denote the class number of the field $\mathbb{Q}(\sqrt{-5l})$.  The proof will be given in the paper \cite{mor3}.  \bigskip

\noindent {\bf Theorem 1.4.} {\it Let $l > 5$ be a prime.}  \begin{enumerate}[a)]
\item {\it Assume $l \equiv 2,3$ mod $5$.  The number of irreducible quartics of the form $f(x) = x^4+ax^3+(11a+2)x^2-ax+1$ dividing $\hat H_{5,l}(x)$ over $\mathbb{F}_l$ is: }
\begin{enumerate}[i)]
\item {\it $\frac{1}{4}h(-5l)$, if $l \equiv 1$ mod $4$;}
\item {\it $\frac{1}{2}h(-5l)-1$, if $l \equiv 3$ mod $8$;}
\item {\it $h(-5l)-1$, if $l \equiv 7$ mod $8$.}
\end{enumerate}
\item {\it If $l \equiv 4$ mod $5$, the number of irreducible quadratic factors $k(x)=x^2+rx+s$ dividing $\hat H_{5,l}(x)$ over $\mathbb{F}_l$, with $r = \varepsilon ^5 (s-1)$ or $r = \bar \varepsilon^5 (s-1)$ in $\mathbb{F}_l$ is:}
\begin{enumerate}[i)]
\item {\it $\frac{1}{2}h(-5l)$, if $l \equiv 1$ mod $4$;}
\item {\it $h(-5l)-3$, if $l \equiv 3$ mod $8$;}
\item {\it $2h(-5l)-3$, if $l \equiv 7$ mod $8$.}
\end{enumerate}
\item {\it If $l \equiv 1$ mod $5$, the number of irreducible quadratic factors $k(x)=x^2+rx+s$ dividing $\hat H_{5,l}(x)$ over $\mathbb{F}_l$, with $r = \varepsilon ^5 (s-1)$ or $r = \bar \varepsilon^5 (s-1)$ in $\mathbb{F}_l$ is:}
\begin{enumerate}[i)]
\item {\it $\frac{1}{2}h(-5l)$, if $l \equiv 1$ mod $4$;}
\item {\it $h(-5l)-1$, if $l \equiv 3$ mod $8$;}
\item {\it $2h(-5l)-1$, if $l \equiv 7$ mod $8$.}
\end{enumerate}
\end{enumerate}
\bigskip

The above theorems are analogues of theorems that have been proven in previous papers for the Hasse invariants of other normal forms for elliptic curves.  Theorems 1.1 and 1.3 are the analogues for the prime $5$ and the Tate normal form $E_5$ of Theorem 1 in \cite{brm}, which gives the counts of linear and certain quadratic factors of the Hasse invariants of the Legendre normal form, the Jacobi normal form, and the Deuring normal form in characteristic $p$, in terms of the class number of the field $\mathbb{Q}(\sqrt{-p})$.  Theorem 1.4 is the analogue of Theorem 1.1 of \cite{mor4}, which gives the count of binomial quadratic factors of the Legendre polynomial $P_{(p-e)/4}(x)$ (mod $p$), where $p$ is an odd prime and $e=1,3$, in terms of the class number of $\mathbb{Q}(\sqrt{-2p})$.  This case is related to the Hasse invariant of the Tate normal form $E_4$ for a point of order $4$.  Theorem 1.4 is also the analogue for the prime $5$ of Theorem 1.1 of \cite{mor5}, which gives the count of binomial quadratic factors of the Legendre polynomial $P_{(p-e)/3}(x)$ (mod $p$) ($e=1,2$), in terms of the class number of the field $\mathbb{Q}(\sqrt{-3p})$; and in this case the proof involves the Hasse invariant of the Deuring normal form $E_3: Y^2+\alpha XY+Y=X^3$. \medskip

The results in \cite{brm}, \cite{mor4}, and \cite{mor5} just referred to have recently been used by Nakaya \cite{na} to count the number of linear factors (mod $p$) of the supersingular polynomial $ss_p^{(N*)}(x)$ corresponding to the Fricke group $\Gamma_0^*(N)$, for $N=2$ and $N=3$ (see \cite{sa1}, \cite{sa2}, \cite{na}).  See Section 6 below for the precise definitions.  In his paper \cite{na} Nakaya also states many interesting conjectures regarding $ss_p^{(N*)}(X)$, including a formula for the number of linear factors of $ss_p^{(N*)}(x)$ for $N \ge 5$ and $N \in \mathfrak{S}-\{2\}$ in his Conjecture 5, where $\mathfrak{S}$ is the set given immediately below.  Using Theorems 1.1 and 1.3, I show in Section 7 that Theorem 1.4, stated above, implies Nakaya's Conjecture 5 for $N=5$.  The proof of Theorem 1.4 appears in a sequel to this paper (see \cite{mor3}), which establishes Nakaya's conjecture for this case. \medskip

In addition, Nakaya \cite{na} also asks the interesting question, whether the roots $j_N^*$ of $ss_p^{(N*)}(X)$ always lie in $\mathbb{F}_{p^2}$, when $N$ is a prime in the set
$$\mathfrak{S} = \{2, 3, 5, 7, 11, 13, 17, 19, 23, 29, 31, 41, 47, 59, 71\}.$$
(See \cite[Question 1]{na}.)  The set $\mathfrak{S}$ consists of the primes $N=q$ for which all supersingular invariants of elliptic curves in characteristic $q$ lie in the prime field $\mathbb{F}_q$.  In Section 6 I give a partial answer to Nakaya's question by showing that the  roots of $ss_p^{(N*)}(X)$ lie in $\mathbb{F}_{p^2}$ for $N \in \{2,3,5,7\}$.  \medskip

For $N=5$, this follows from Theorem 6.2 of \cite{mor}, which concerns the form of the irreducible factors of $\hat H_{5,l}(x)$ modulo $l$.  In order to prove this for $N=7$, I consider the Hasse invariant $\hat H_{7,l}(x)$ of the Tate normal form $E_7$ for a point of order $7$.  In Theorems 5.2 and 5.3 I determine the degrees of the irreducible factors of $\hat H_{7,l}(x)$ modulo $l$, as well as the particular forms of cubic factors and sextic factors which can appear.   The proofs of Theorems 5.2 and 5.3 make use of meromorphisms of elliptic function fields, as in \cite{mor}.  Conjecture 1 in Section 5 is the analogue of Theorem 1.1 for the polynomial $\hat H_{7,l}(x)$.  \medskip

The paper is laid out as follows. In Section 2 I handle the cases $j=0, 1728$ of Theorem 1.2, and collect the preliminary results on the icosahedral group that will be needed for the proofs of Theorems 1.1-1.3.  Section 3 contains the proof of Theorem 1.2.  In Section 4 I determine the factorization of the polynomial $p_d(x)$ (mod $l$), for $d=l$ or $d=4l$, in case $l \equiv 1, 4$ (mod $5$), and I prove Theorem 1.3 (see Theorem 4.3).  Section 5  discusses the Hasse invariant of the Tate normal form $E_7$, as well as Conjecture 1 on the number of linear and cubic factors of $\hat H_{7,l}(X)$ and Conjecture 2 on the number of linear factors of $\hat H_{p,l}(X)$ over $\mathbb{F}_l$, for a general prime $p$.  Theorems 5.2 and 5.3 are used in Section 6 to give a partial affirmative answer to Nakaya's question mentioned above.  Finally, Theorem 7.1 is proved in Section 7, showing that Theorem 1.4 implies Nakaya's Conjecture 5 for $ss_p^{(5*)}(X)$ in \cite{na}.

\section{Preliminaries.}

\subsection{The $j$-invariants $j=0, 1728$.}

We start by proving the assertion of Theorem 1.2 for the two polynomials $f_0(x), f_{1728}(x)$.  From the representations
$$f_0(x)=(x^2+6x-1)^2-20x^2, \ \ f_{1728}(x)=(x^2+9x-1)^2-5x^2,$$
it is not difficult to see that the splitting fields of these two polynomials over $\mathbb{Q}$ are
$$L_1=\mathbb{Q}\left(\sqrt{15-6\sqrt{5}}\right), \ \ L_2=\mathbb{Q}\left(\sqrt{10-2\sqrt{5}}\right),$$
respectively.  These square-roots can be represented in terms of $15$th and $20$th roots of unity:
$$\sqrt{15-6\sqrt{5}}=1-4\zeta_{15}+2\zeta_{15}^2+2\zeta_{15}^5-2\zeta_{15}^6-2\zeta_{15}^7, \ \ \zeta_{15}=e^{\frac{2\pi i}{15}},$$
$$\sqrt{10-2\sqrt{5}}=2\zeta_{20}^3-2\zeta_{20}^7, \ \ \zeta_{20}=e^{\frac{2\pi i}{20}}.$$
Since $L_1$ and $L_2$ are totally real, it follows that
$$L_1=\mathbb{Q}(\zeta_{15}+\zeta_{15}^{-1}), \ \ L_2=\mathbb{Q}(\zeta_{20}+\zeta_{20}^{-1})$$
are the real subfields of $\mathbb{Q}(\zeta_{15})$ and $\mathbb{Q}(\zeta_{20})$, respectively.  \medskip

This gives the following proposition, in which we consider the factorization of the polynomials $f_0, f_{1728}$ modulo primes which are $\pm 1$ (mod $5$).  \bigskip

\noindent {\bf Proposition 2.1.} \begin{enumerate}[a)]
\item {\it If $l \equiv 4$ (mod $5$) is a prime for which $\left(\frac{-3}{l}\right)=-1$, so that $j=0$ is supersingular in characteristic $l$, then $f_0(x)$ splits into distinct linear factors (mod $l$).} 
\item {\it If $l \equiv 4$ (mod $5$) is a prime for which $\left(\frac{-4}{l}\right)=-1$, so that $j=1728$ is supersingular in characteristic $l$, then $f_{1728}(x)$ splits into distinct linear factors (mod $l$).} 
\item {\it If $l \equiv 1$ (mod $5$) is a prime for which $\left(\frac{-3}{l}\right)=-1$, then $f_0(x)$ splits into two irreducible quadratics of the form $q(x) =x^2+ax-1$ (mod $l$).} 
\item{\it If $l \equiv 1$ (mod $5$) is a prime for which $\left(\frac{-4}{l}\right)=-1$, then $f_{1728}(x)$ splits into two irreducible quadratics of the form $q(x) =x^2+ax-1$ (mod $l$).} 
\end{enumerate} \medskip

\noindent {\it Proof.} a) If $l \equiv 4$ (mod $5$) and $\left(\frac{-3}{l}\right)=-1$, then $l \equiv -1$ (mod $15$), so that $l$ splits completely in the field $L_1$.  Note that $l$ does not divide $\textrm{disc}(f_0(x))=
2^{12}3^2 5^3$.  A similar statement holds for part b), since $l \equiv -1$ (mod $20$) and $\textrm{disc}(f_{1728}(x))=2^8 3^4 5^3$.  \medskip

To prove c), note that $l \equiv 1$ (mod $5$) and $l \equiv 2$ (mod $3$) implies $l \equiv 11$ (mod $15$).  Since the order of $11$ with respect to the subgroup $\{\pm1\}$ modulo $15$ is $2$, $f_0(x)$ splits into two irreducible quadratics (mod $l$), corresponding to the two prime divisors of $l$ of degree $2$ in $L_1$.  That these quadratics have the given form follows from the factorization of $f_0(x)$ over the subfield $\mathbb{Q}(\sqrt{5})$:
$$f_0(x) = (x^2+(6+2\sqrt{5})x-1)(x^2+(6-2\sqrt{5})x-1).$$
The proof of d) is similar, using that $l \equiv 11$ (mod $20$) and
$$f_{1728}(x) = (x^2+(9+\sqrt{5})x-1)(x^2+(9-\sqrt{5})x-1).$$
$\square$  \bigskip

From \cite{mor} we have the following alternate representation of the Hasse invariant for $E_5$:
\begin{align*}
\hat H_{5,l}(b)= & g_0(b)^r (b^2+1)^s g_{1728}(b)^sb^{n_l}(1-11b-b^2)^{5n_l}\\
& \times J_l\left(\frac{(b^4-228b^3+494b^2+228b+1)^3}{b(1-11b-b^2)^5}\right),
\end{align*}
where
\begin{align*}
g_0(b)= & b^4-228b^3+494b^2+228b+1,\\
g_{1728}(b)= & b^4+522b^3-10006b^2-522b+1.
\end{align*}
Proposition 2.1 also holds with $g_0(x)$ in place of $f_0(x)$ and $g_{1728}(x)$ in place of $f_{1728}(x)$.  To see this, consider the equations
$$g_0(x)=(x^2-114x-1)^2-2^2 5^5x^2, \ \ \ g_{1728}(x) = (x^2+261x-1)^2-5^7x^2.$$
These formulas imply that the splitting fields of these two polynomials are $\mathbb{Q}(\sqrt{255+114\sqrt{5}})$ and $\mathbb{Q}(\sqrt{650-290\sqrt{5}})$, respectively.  Since
\begin{align*}
255+114\sqrt{5} & =\left(\frac{1+\sqrt{5}}{2}\right)^{12}(15-6\sqrt{5})\\
650-290\sqrt{5} & = (\sqrt{5})^2\left(\frac{1-\sqrt{5}}{2}\right)^6(10-2\sqrt{5}),
\end{align*}
it follows from Kummer theory that these splitting fields coincide with the fields $L_1$ and $L_2$ above.  Because $\textrm{disc}(g_0(x))=2^{12} 3^2 5^{15}$ and $\textrm{disc}(g_{1728}(x))=2^8 3^4 5^{21}$, the polynomials $g_0(x), g_{1728}(x)$ factor over $\mathbb{F}_l$ in the same way that $f_0(x), f_{1728}(x)$ do.  In particular, we have the factorizations
\begin{align*}
g_0(x) & = (x^2-(114+50\sqrt{5})x-1)(x^2-(114-50\sqrt{5})x-1),\\
g_{1728}(x) & = (x^2 + (261+125\sqrt{5})x-1)(x^2 + (261-125\sqrt{5})x-1).
\end{align*}

\subsection{The group $G_{60}$.}

The above representation of the Hasse invariant will be useful, because the rational function
$$j_5(z)=\frac{(1+228z^5+494z^{10}-228z^{15}+z^{20})^3}{z^5(1-11z^5-z^{10})^5} \eqno{(2.1)}$$
is invariant under a group $G_{60}$ of linear fractional substitutions:
$$G_{60}= \langle S,T \rangle, \ \ S(z)=\zeta z, \ \ T(z)=\frac{-(1+\sqrt{5})z+2}{2z+1+\sqrt{5}}, \ \ (\zeta = \zeta_5 = e^{2\pi i/5}),$$
which is isomorphic to the icosahedral group $A_5$.  The coefficients of the maps in $G_{60}$ are in the field $\mathbb{Q}(\zeta_5)$ and $S$ and $T$ have orders 5 and 2, respectively (see \cite{hup}).  The transformation $T$, together with the transformation $U(z)=\frac{-1}{z}$, given in terms of $S$ and $T$ by $U=T \cdot S^2 \cdot T \cdot S^3 \cdot T \cdot S^2$ (see \cite[II, pp. 42-43]{fr1}), generates the subgroup
$$H=\{1, T, U, TU\}$$
of $G_{60}$, where $TU(z)=UT(z)=-1/T(z)=T_2(z)$, and
$$T_2(z)=\frac{-(1-\sqrt{5})z+2}{2z+1-\sqrt{5}}.$$
Thus, $U=TT_2=T_2T$.  The normalizer of $H$ in $G_{60}$ is $N=\langle A, H \rangle \cong A_4$, where $A=STS^{-2}$ is the map
$$A(z)=\zeta^3 \frac{(1+\zeta)z+1}{z-1-\zeta^4}$$
of order $3$, and $ATA^{-1}=U, AUA^{-1}=T_2$.  Also, $A^\sigma=A^{-1}U$ is the conjugate map
$$A^\sigma(z)=\zeta \frac{(1+\zeta^2)z+1}{z-1-\zeta^3},$$
obtained by applying the automorphism $\sigma: \zeta \rightarrow \zeta^2$ to the coefficients.  Thus,
$$A^{2}=A^{-1}=A^{\sigma}U.$$

Let $F(x,j)$ be the polynomial in characteristic $0$ given by
$$F(x,j)=(x^4+12x^3+14x^2-12x+1)^3-jx^5(1-11x-x^2).$$
Then $F(x,j)$ is a factor of
$$F_d(x)=x^{5\textsf{h}(-d)}(1-11x-x^2)^{\textsf{h}(-d)}H_{-d}\left(\frac{(x^4+12x^3+14x^2-12x+1)^3}{x^5(1-11x-x^2)}\right) \eqno{(2.2)}$$
whenever $j$ is a root of the class equation $H_{-d}(x)$ corresponding to the discriminant $-d=-l$ or $-4l$.  From \cite{mor} we recall that the polynomials
$$f_0(x)=x^4+12x^3+14x^2-12x+1, \ \ g_0(x)=x^4-228x^3+494x^2+228x+1$$
and
$$\Delta(x)=x^5(1-11x-x^2), \ \ \Delta_5(x)=x(1-11x-x^2)^5$$
have the property that
$$(\varepsilon^5 x+1)^4 f_0(\tau(x))=5\varepsilon^{10}g_0(x), \ \ (\varepsilon^5 x+1)^{12} \Delta(\tau(x))=125\varepsilon^{30} \Delta_5(x),$$
where $\tau(x)=\frac{-x+\varepsilon^5}{\varepsilon^5x+1}$ and $\varepsilon = \frac{-1+\sqrt{5}}{2}$.  It follows that
$$(\varepsilon^5 x+1)^{12} F(\tau(x),j)=5^3 \varepsilon^{30}G(x,j),$$
where
$$G(x,j)=(x^4-228x^3+494x^2+228x+1)^3-jx(1-11x-x^2)^5.\eqno{(2.3)}$$
The latter polynomial is a factor of the polynomial $G_d(x)$ defined by
$$G_d(x)=x^{\textsf{h}(-d)}(1-11x-x^2)^{5\textsf{h}(-d)}H_{-d}\left(\frac{(x^4-228x^3+494x^2+228x+1)^3}{x(1-11x-x^2)^5}\right),$$
for the same values of $j$.  It follows (when $l \equiv 4$ mod $5$) that there is a 1-1 correspondence between linear factors of $F(x,j) \equiv f_j(x)$ over $\mathbb{F}_l$ and linear factors of $G(x,j)$ over $\mathbb{F}_l$, since $\sqrt{5} \in \mathbb{F}_l$ and the determinant of the linear fractional map $\tau(x)$ is
$$\textrm{det}(\tau(x))=-1-\varepsilon^{10}=\frac{-125+55\sqrt{5}}{2}=-5\sqrt{5}\varepsilon^5.$$

The polynomial $G(x^5,j)$ satisfies the transformation formulas
$$G(S(x)^5,j)=G(x^5,j),$$
$$ (2x+1+\sqrt{5})^{60}G(T(x)^5,j)=2^{60}\left(\frac{5+\sqrt{5}}{2} \right)^{30} G(x^5,j).$$
Thus, for a given value of $j$, the roots of $G(x^5,j)=0$ are invariant (as a set) under the group $G_{60}$, and this is also true over $\mathbb{F}_l$, since $l > 5$.  Moreover, the number of linear factors of $G(x^5,j)$ (counted with multiplicity) is the same as the number of linear factors of $G(x,j)$ (counted with multiplicity), and therefore of $F(x,j)$, since $x \rightarrow x^5$ is a bijection on $\mathbb{F}_l$ (for $l \equiv 4$ mod $5$) and $\zeta \in \mathbb{F}_{l^2}-\mathbb{F}_l$.  Also note that the factor $2x+1+\sqrt{5}$ is not $0$ in $\mathbb{F}_l$, for any root $x$ of the congruence $G(x^5,j) \equiv 0$ (mod $l$), because $(x+\frac{1+\sqrt{5}}{2})(x+\frac{1-\sqrt{5}}{2})=x^2+x-1$ is a factor of $\Delta_5(x^5)$; and this is not $0$ (mod $l$), since $j$ is a supersingular invariant of $E_5$ in characteristic $l$.  \medskip

To complete the proof of Theorem 1.2, we just need to show that $G(x^5,j)$ has exactly four distinct linear factors (mod $l$), when $j \in \mathbb{F}_l$ is a supersingular $j$-invariant different from $0$ and $1728$.  Finally, note that $G(x^5,j)$ has no multiple factors (mod $l$) for these values of $j$, since
$$\textrm{disc}_x(G(x,j))=5^{145}j^8(j-1728)^6.$$

\section{Proof of Theorem 1.2.}

We shall make use of the polynomial $p_d(x)$ from \cite{mor2}, where $d=l$ or $d=4l$ and $l \equiv 4$ (mod $5$) is a prime.  From that paper we know that $p_d(x)$ is a factor of both polynomials $F_d(x)$ and $G_d(x)$ (see \cite[p. 19]{mor2}) and is stabilized by the group $H$, meaning that its roots are invariant under the mappings in $H$:
$$\textrm{Stab}_{G_{60}}(p_d(x)) = H =  \{1, T, U, TU\}.$$
Here the action of $G_{60}$ on a polynomial $p(x)$ is defined by
$$\sigma(x) = \frac{ax+b}{cx+d} \ \longrightarrow \ p^\sigma(x) = (cx+d)^{deg(p)}p(\sigma(x)).$$
It can be checked that
$$\left(z+\frac{1+\sqrt{5}}{2}\right)^{4\textsf{h}(-d)}p_d(T(z))= \left(\frac{5+\sqrt{5}}{2} \right)^{2\textsf{h}(-d)} p_d(z).$$
(See \cite[p. 30]{mor2}.) \medskip

Let $K=\mathbb{Q}(\sqrt{-l})$ and let $R_K$ be its ring of integers, so the ideal $(5)=\wp_5\wp_5'$ splits into a product of the prime ideals $\wp_5, \wp_5'$ in $R_K$.  Let $F_1=\mathbb{Q}(\eta)$, as in \cite{mor2}, where
$$\eta=r(w/5), \ \ \textrm{with} \ \ w=\frac{v+\sqrt{-d}}{2}, \ \ \wp_5^2 \mid w, \ \ d=l \ \textrm{or} \ 4l,\eqno{(3.1)}$$
and $r(\tau)$ is the Rogers-Ramanujan continued fraction.  (See \cite{anb}, \cite{ber}, \cite{du}. As in \cite{mor2}, we also need $(N(w),2)=1$ in case $d=4l$ and $l \equiv 3$ mod $4$.)  Then $\eta$ is a root of the irreducible polynomial $p_d(x)$ and $F_1$ is an abelian extension of the quadratic field $K$.  Let $\mathcal{O}=\textsf{R}_{-d}$ be the order of discriminant $-d=-l$ or $-d=-4l$ in $K$ (the former only for $l \equiv 3$ (mod $4$)).  It is shown in \cite{mor2} that
$$F_1=\Sigma_{\wp_5'}, \ \ \textrm{if} \ \ d=l, \ \ \textrm{or} \ d=4l \ \textrm{and} \ l \equiv 1 \ (\textrm{mod} \ 4),$$
where $\Sigma_{\wp_5'}$ is the ray class field of conductor $\mathfrak{f}=\wp_5'$ over $K$.  On the other hand,
$$F_1=\Sigma_{\wp_5'}\Omega_2, \ \  \textrm{if} \ \ d= 4l \ \textrm{and} \ l \equiv 3 \ (\textrm{mod} \ 4),$$
where $\Omega_2$ is the ring class field of conductor $2$ over $K$. \medskip

The value $\xi=\zeta^j r(-1/w)$, for some $j \not \equiv 0$ (mod $5$), also generates a class field over $K$, namely the conjugate field
$$F_2 = \mathbb{Q}(\xi) = \Sigma_{\wp_5} \ \ \textrm{or} \ \ \Sigma_{\wp_5}\Omega_2$$
of $F_1$ in the two cases described above.  It is shown in \cite{mor2} that $F_1 \neq F_2$ are disjoint quadratic extensions of the Hilbert class field $\Sigma_1$ in the first case; and of $\Omega_2$ in the second.  Moreover, the unit $\xi$ is also a root of $p_d(x)$, and is given by the formula
$$\xi = \frac{-(1+\sqrt{5})\eta^{\tau_5}+2}{2\eta^{\tau_5}+1+\sqrt{5}}, \ \ \tau_5 =\left(\frac{F_1/K}{\wp_5}\right).$$
Finally, if $d=lf^2$, with $f=1$ or $2$, the Galois group $\textrm{Gal}(F_1/\Omega_f)$ of the quadratic extension $F_1/\Omega_f$ is generated by the automorphism $\sigma=(\eta \rightarrow -1/\eta)$. \smallskip

With these facts in hand we prove the following. \bigskip

\noindent {\bf Theorem 3.1.} {\it If $l$ is a prime with $l \equiv 4$ (mod $5$) and $d=l$ or $4l$, then the polynomial $p_d(x)$ factors (mod $l$) as a product of $2\textsf{h}(-d)$ squares of linear factors.} \medskip

\noindent {\it Proof.} This follows as in \cite{brm} (see Propositions 9 and 10).  This is because $l$ ramifies in $K=\mathbb{Q}(\sqrt{-l})$, so it is the square of a product of prime ideals in $F_1$.  Also, the prime divisor $\mathfrak{l}=(\sqrt{-l})$ of $l$ in $R_K$ is unramified in $F_1$ and $\pm \sqrt{-l}$ is congruent to 1 (mod $\wp_5'$), so $\mathfrak{l}$ is in the ray mod $\wp_5'$ and in the principal ring class mod $2$ if $d=2^2l$ and $l \equiv 3$ (mod $4$).  Therefore, $\mathfrak{l}$ splits into a product of primes of degree 1 in $F_1$ (and ramification index $e=2$ over $\mathbb{Q}$).  Since a root $\eta$ of $p_d(x)$ generates $F_1$ over $\mathbb{Q}$, $p_d(x)$ factors over the $l$-adic field $\mathbb{Q}_l$ as a product of irreducible quadratics, each of which belongs to a ramified extension of $\mathbb{Q}_l$ and therefore has discriminant divisible by $l$.  Thus, $p_d(x) \equiv (x-a_1)^2 \cdots (x-a_{2h})^2$ (mod $l$), where $h=\textsf{h}(-d)$ is the class number of the order $\textsf{R}_{-d}$. $\square$ \bigskip

We will see below that the roots $a_i$ of $p_d(x)$ (mod $l$) in the proof of Theorem 3.1 are distinct. At any rate, $p_d(x)$ is divisible by at most $2\textsf{h}(-d)$ distinct linear factors (mod $l$). \bigskip

\noindent {\bf Proposition 3.2.} {\it If $j \neq 0, 1728$ is a supersingular $j$-invariant (mod $l$), then the congruence $G(x^5, j) \equiv 0$ (mod $l$) has at least $4$ distinct solutions in $\mathbb{F}_l$, which are also solutions of $p_d(x) \equiv 0$ (mod $l$), for a unique value of $d=l$ or $4l$, depending on $j$.} \medskip

\noindent {\it Proof.} Let $j \neq 0, 1728$ be a supersingular $j$-invariant (mod $l$).  Choose a root $\mathfrak{j}$ of $H_{-d}(x)=0$ in characteristic $0$, with $d=l$ or $4l$, such that $\mathfrak{j}$ reduces to $j$ modulo some prime divisor $\mathfrak{q}$ of $l$ in the field $F_1$.  This is possible because of the congruences (1.1) and (1.2) for $H_{-d}(x)$ given in the introduction and Deuring's lifting theorem.  Let $\Sigma = K(\mathfrak{j})$, so this field is either the Hilbert class field of $K$, for $d=l$, or $d=4l$ and $l \equiv 1$ (mod $4$); or $\Sigma=\Omega_2$ is the ring class field over $K$ with conductor $f=2$, when $d=4l$ and $l \equiv 3$ (mod $4$).  Then $j \equiv \mathfrak{j}$ also modulo $\mathfrak{p}$, where $\mathfrak{p}$ is the prime divisor of $\Sigma$ which $\mathfrak{q}$ divides. \medskip

From \cite[eq. (25)]{mor2} we know that if $w$ is chosen as in (3.1), then
$$j(w/5) = \frac{(\eta^{20}-228\eta^{15}+494\eta^{10}+228\eta^5+1)^3}{\eta^5(1-11\eta^5-\eta^{10})^5}.$$
(Also see \cite[ eq. (2.5)]{du}.)  Since the root $\mathfrak{j}$ of $H_{-d}(x)=0$ is a conjugate of $j(w/5)$, some conjugate of $\eta$ over $K$ is a root of $G(x^5,\mathfrak{j})=0$.  Thus, the minimal polynomial $\mu_1(x)$ of some $\eta^\sigma$ over $\Sigma$ (with $\sigma \in \textrm{Gal}(F_1/K)$) divides $G(x^5,\mathfrak{j})$.  Since $\mathfrak{j} \equiv j$ (mod $\mathfrak{q}$), then
$$G(\eta^{5\sigma},j) \equiv G(\eta^{5\sigma},\mathfrak{j}) \equiv 0 \ (\textrm{mod} \ \mathfrak{q});$$
and the fact that $\textrm{deg}_\mathbb{Q}(\mathfrak{q})=1$ means that $G(x^5, j) \equiv 0$ (mod $\mathfrak{p}$) has at least two distinct solutions $a_1, a_2=-1/a_1$ (mod $\mathfrak{p}$) in $\mathbb{F}_l$ (recalling that roots of $G(x^5,j)$ are invariant under the mapping $U(z)$).  From \cite[p. 1193]{mor2} we also have the equation
$$j(w)= \frac{(\xi^{20}-228\xi^{15}+494\xi^{10}+228\xi^5+1)^3}{\xi^5(1-11\xi^5-\xi^{10})^5}.$$
Thus, the minimal polynomial $\mu_2(x)$ over $\Sigma$ of some conjugate of $\xi$ over $K$ also divides $G(x^5,\mathfrak{j})$.  Now, $\mu_1(x)$ and $\mu_2(x)$ are certainly relatively prime to each other, since their respective roots generate distinct extensions $F_1 \neq F_2$ over $\Sigma$.  Thus, $\mu_1(x)\mu_2(x) \mid G(x^5,\mathfrak{j})$ in $\Sigma[x]$.  Since the congruence $G(x^5,j) \equiv 0$ does not have any multiple roots (mod $l$), this shows that $G(x^5,j)$ has at least four distinct linear factors modulo $\mathfrak{p}$ and therefore also modulo $l$.  We know that $\mu_1(x)\mu_2(x) \mid p_d(x)$, so these linear factors are also factors of $p_d(x)$ (mod $l$). $\square$ \medskip

\noindent {\bf Remark.}  The same arguments as in the above proof show that $p_d(x) \equiv 0$ (mod $l$) has at least one solution in $\mathbb{F}_l$ corresponding to each of the values $j=0$ and $j=1728$.  \medskip

By the discussion in Section 2, the roots of $G(x^5,j)=0$ are invariant under the group $G_{60}$, so $G(x^5,\mathfrak{j})$ splits into linear factors in the field $F=\mathbb{Q}(\eta, \zeta)$.  If $\overline{\mathfrak{q}}$ is a prime divisor of $\mathfrak{q}$ in the field $F$, all the roots of the congruence $G(x^5, j) \equiv 0$ (mod $\overline{\mathfrak{q}}$) are images of the root $a_1$, say, by elements of $\bar G_{60} = G_{60}$ (mod $\overline{\mathfrak{q}}$).  Note that $G_{60} \cong A_5$ is a simple group, so the kernel of the map taking $G_{60}$ to its reduction modulo $\overline{\mathfrak{q}}$ in $F$ is trivial.  Four of the images $M(a_1)$ lie in $\mathbb{F}_l$, since the elements of the group $H$ have coefficients in $\mathbb{Q}(\sqrt{5})$.  We now prove the following.  \bigskip

\noindent {\bf Proposition 3.3.} {\it If $a_1 \in \mathbb{F}_l$ is a root of the congruence $G(x^5, j) \equiv 0$ (mod $l$), with $j \not \equiv 0, 1728$ (mod $l$), then $M(a_1) \in \mathbb{F}_{l^2}-\mathbb{F}_l$, for $M \in \bar G_{60}- \bar H$.} \medskip

Before proceeding with the proof, we recall (2.1) and record the factorizations
\begin{align*}
j_5(x)=& \frac{(x^4-3x^3-x^2+3x+1)^3h_1(x)^3h_2(x)^3}{\Delta_5(x^5)},\\
j_5(x)-1728=& \frac{(x^2+1)^2(x^4+2x^3-6x^2-2x+1)^2h_3(x)^2h_4(x)^2h_5(x)^2}{\Delta_5(x^5)},
\end{align*}
where the $h_i(x)$ are irreducible polynomials of degree $8$ (for $1 \le i \le 5$) and
\begin{align*}
\Delta_5(x^5)&=-x^5(x^2+x-1)^5f(x)^5 \times x^{20}f(-1/x)^5,\\
f(x)&=x^4+2x^3+4x^2+3x+1.
\end{align*}
Also note the resultant
$$Res_x(x^4-228x^3+494x^2+228x+1, 1-11x-x^2)=5^{10};$$
this shows that $\Delta_5(x^5)$ is never $0$ modulo $l$ for a root of $G(x^5,j) \equiv 0$ (mod $l$). \medskip

\noindent {\it Proof of Proposition 3.3.} To prove the proposition, it suffices to show $M(z) \not \in \mathbb{F}_l$ for any root $z \in \mathbb{F}_l$ of $G(z^5,j) \equiv 0$ (mod $l$), for a system $\{M\}$ of left coset representatives of $\bar H$ in $\bar G_{60}$.  It is not hard to see that such a system is given by
$$M_{ik}=S^k A^i, \ \ 0 \le i \le 2, \ 0 \le k \le 4,$$
where $A, S$ are the linear fractional maps from Section 2.2.  First, if $i=0$ and $k \neq 0$, $S^{k}(z)=\zeta^{k}z$ cannot lie in $\mathbb{F}_l$ if $z \neq 0$ does, since $\zeta$ is quadratic over $\mathbb{F}_l$ (this because $l \equiv 4$ mod $5$). \medskip

If $x=A(z)$, then $x \in \mathbb{F}_l$ if and only if
$$0 \equiv x-x^l \equiv A(z)-A(z)^l \equiv \frac{(2\zeta^2 +2\zeta +1)(z^2+2(\zeta^2+\zeta^3+1)z-1)}{(-z+1+\zeta)(z+\zeta+\zeta^2+\zeta^3)} \ (\textrm{mod} \ \overline{\mathfrak{q}}).$$
For this, note that $A(z)^l \equiv \zeta^2 \frac{(1+\zeta^4)z+1}{z-1-\zeta}$ (mod $\overline{\mathfrak{q}}$), since $\zeta^l = \zeta^4$.  However, $N_{\mathbb{Q}(\zeta)/\mathbb{Q}}(2\zeta^2 +2\zeta +1)=5$, and
$$(z^2+2(\zeta^2+\zeta^3+1)z-1)(z^2+2(\zeta^4+\zeta+1)z-1)=z^4+2z^3-6z^2-2z+1,$$
which cannot be $0$, since the second of the above factorizations would imply that $z$ satisfies $j_5(z)=1728$ and therefore $G(z^5,1728) \equiv 0$, contradicting $j \not \equiv 1728$ (mod $l$).  Hence, $x^l \not \equiv x$ (mod $\overline{\mathfrak{q}}$), giving that $x =A(z) \not \in \mathbb{F}_l$.  \medskip

The arguments for the other coset representatives $M_{ik} \neq M_{00}$ are similar: modulo $\overline{\mathfrak{q}}$ we have the congruences
\begin{align*}
SA(z)-SA(z)^l &\equiv \frac{(\zeta^3+\zeta^2+2\zeta+1)(z-\zeta^2-\zeta^3)(z+\zeta^3+\zeta^2+1)}{(-z+1+\zeta)(z+\zeta+\zeta^2+\zeta^3)}\\
S^2A(z)-S^2A(z)^l &\equiv \frac{-(\zeta^3+\zeta^2+2\zeta+1)(z^2+(\zeta^2+\zeta^3-1)z-1))}{(-z+1+\zeta)(z+\zeta+\zeta^2+\zeta^3)}\\
S^3A(z)-S^3A(z)^l &\equiv \frac{-(2\zeta^2 +2\zeta +1)(z^2-(\zeta^2+\zeta^3+2)z-1)}{(-z+1+\zeta)(z+\zeta+\zeta^2+\zeta^3)}\\
S^4A(z)-S^4A(z)^l &\equiv \frac{(\zeta^3-3\zeta^2-2\zeta-1)z}{(-z+1+\zeta)(z+\zeta+\zeta^2+\zeta^3)}.
\end{align*}
The norms of the constant multipliers in the numerators are, respectively, $5, 5, 5, 5^3$; and the $z$-factors in the numerators divide the respective polynomials
$$(z^2+z-1)^2, \ \ z^4-3z^3-z^2+3z+1, \ \ z^4-3z^3-z^2+3z+1, \ \ z.$$
Since $z(z^2+z-1) \mid \Delta_5(z^5)$ and $(z^4-3z^3-z^2+3z+1) \mid j_5(z)$, none of them can be $0$ for a root of $G(z^5,j) \equiv 0$.  Hence $x = S^kA(z) \not \in \mathbb{F}_l$ for $0 \le k \le 4$.
 \medskip

Taking $\displaystyle A^2(z)=-\zeta \frac{(1+\zeta+\zeta^2)z-\zeta^2}{z+1+\zeta+\zeta^2}$, we also have the following congruences modulo $\overline{\mathfrak{q}}$:
\begin{align*}
A^2(z)-A^2(z)^l &\equiv \frac{(\zeta^3+\zeta^2+2\zeta+1)(z^2-2(\zeta^2+\zeta^3)z-1)}{(\zeta^2+\zeta-z)(z+1+\zeta+\zeta^2)}\\
SA^2(z)-SA^2(z)^l &\equiv \frac{-(\zeta^3+\zeta^2+2\zeta+1)(z^2+(\zeta^2+\zeta^3-1)z-1)}{(\zeta^2+\zeta-z)(z+1+\zeta+\zeta^2)}\\
S^2A^2(z)-S^2A^2(z)^l &\equiv \frac{(2\zeta+2\zeta^2+1)(\zeta^2+\zeta^3-z)(z+1+\zeta^2+\zeta^3)}{(\zeta^2+\zeta-z)(z+1+\zeta+\zeta^2)}\\
S^3A^2(z)-S^3A^2(z)^l &\equiv \frac{-(2\zeta^3+4\zeta^2+6\zeta+3)z}{(\zeta^2+\zeta-z)(z+1+\zeta+\zeta^2)}\\
S^4A^2(z)-S^4A^2(z)^l &\equiv \frac{(2\zeta+2\zeta^2+1)(z^2-(\zeta^2+\zeta^3+2)z-1)}{(\zeta^2+\zeta-z)(z+1+\zeta+\zeta^2)}.
\end{align*}
As before, the norms of the constant multipliers are nonzero (mod $l$) (the fourth norm in the sequence is $5^3$), and the $z$-factors divide the respective polynomials
$$z^4+2z^3-6z^2-2z+1, \ \ z^4-3z^3-z^2+3z+1, \ \ (z^2+z-1)^2, \ \ z, \ \ z^4-3z^3-z^2+3z+1.$$
This gives that $x = S^k A^2(z) \not \in \mathbb{F}_l$ for $0 \le k \le 4$. \medskip

It follows that $M(a_1) \not \in \mathbb{F}_l$ for any $M \not \in \bar H$. $\square$ \bigskip

Propositions 3.2 and 3.3 now yield that $G(x^5,j) \equiv 0$ (mod $l$) has exactly $4$ distinct roots in $\mathbb{F}_l$.  Together with Proposition 2.1, this completes the proof of Theorem 1.2.

\section{A discriminant theorem for $p_d(x)$.}

Next, we derive the following theorem from the considerations in Section 3.  \bigskip

\noindent {\bf Theorem 4.1.} {\it If $d=l$ or $4l$, where $l \equiv 4$ (mod $5$) is prime, the polynomial $p_d(x)$ factors modulo $l$ as
$$p_d(x) \equiv \prod_{i=1}^{2h}{(x-a_i)^2} \ \ (\textrm{mod} \ l),$$
where the $a_i$ are distinct elements of $\mathbb{F}_l$ and $h=\textsf{h}(-d)$ is the class number of the order $\textsf{R}_{-d}$ of discriminant $-d$ in $K=\mathbb{Q}(\sqrt{-l})$.} \medskip

\noindent {\it Proof.} From the results of Section 3, we know that $p_d(x)$ has at most $2h$ distinct roots (mod $l$), and that these roots lie in $\mathbb{F}_l$.  We also know that for every supersingular $j$-invariant in $\mathbb{F}_l$, different from $j=0$ or $j=1728$, which is a root of the congruence $H_{-d}(x) \equiv 0$ (mod $l$), there are $4$ distinct roots of $p_d(x) \equiv 0$ (mod $l$) in $\mathbb{F}_l$.  We now prove the same thing for $j=0$ and $d=4l$.  (Note that $j=0$ is never a root of the class equation $H_{-l}(x)$ modulo $l$.  See the proof of Proposition 11 in \cite{brm}.)  By the remark following the proof of Proposition 3.2, we know that $p_d(x) \equiv 0$ has at least one root $a \in \mathbb{F}_l$ corresponding to $j = 0$, and this is a zero (mod $l$) of
$$G(x^5,0) = g_0(x^5)^3 = (x^{20}-228x^{15}+494x^{10}+228x^5+1)^3.$$
Using the fact that $p_d(x)$ is stabilized by the group $H$, we show that $p_d(x) \equiv 0$ has exactly $4$ distinct roots in $\mathbb{F}_l$ corresponding to $j=0$.  Certainly there are no more than $4$ such roots, since $\textrm{deg}(g_0(x))=4$.  To prove this, consider the fixed points of the mappings in the group $H$.  We compute modulo $l$ that:
\begin{align*}
U(x) \equiv x & \ \ \Rightarrow \ \ x^2 +1 \equiv 0;\\
T(x) \equiv x & \ \ \Rightarrow \ \ x^2+(1+\sqrt{5})x-1 \equiv 0 \ \ \Rightarrow \ x^4+2x^3-6x^2-2x+1 \equiv 0;\\
T_2(x) \equiv x & \ \ \Rightarrow \ \ x^2+(1-\sqrt{5})x-1 \equiv 0 \ \ \Rightarrow \ x^4+2x^3-6x^2-2x+1 \equiv 0.
\end{align*}
However, these factors all divide the numerator of $j_5(x)-1728$, so that any fixed point of $U, T$ or $T_2$ satisfies $j_5(x) \equiv 1728$ (mod $l$).  It follows that there are four distinct images of $a \in \mathbb{F}_l$ which are roots of $p_d(x) \equiv 0$,  and since $\bar H$ stabilizes the roots of $G(x^5,0) \equiv 0$, our claim is proved.  \medskip

By the congruences for $H_{-d}(X)$ in Section 1, there are either $\frac{\textsf{h}(-l)-1}{2}, \frac{\textsf{h}(-4l)-1}{2},$ or $\frac{\textsf{h}(-4l)}{2}$ distinct supersingular $j$-invariants, different from $j=1728$, which are roots of $H_{-d}(x) \equiv 0$ (mod $l$), depending on whether $d=l, 4l$ (with $l \equiv 3$ mod $4$) or $d=4l$ and $l \equiv 1$ modulo $4$.  Thus, $p_d(x) \equiv 0$ has at least
\begin{align*}
4\frac{\textsf{h}(-l)-1}{2} &= 2\textsf{h}(-l)-2, \ (\textrm{when} \ d=l),\\
4\frac{\textsf{h}(-4l)-1}{2} &= 2\textsf{h}(-4l)-2, \ (\textrm{when} \ d=4l, \ l \equiv 3 \ \textrm{mod} \ 4), \ \textrm{or}\\
4\frac{\textsf{h}(-4l)}{2} &=2\textsf{h}(-4l), \ (\textrm{when} \ d=4l, \ l \equiv 1 \ \textrm{mod} \ 4),
\end{align*}
distinct roots corresponding to supersingular $j$-invariants $j \not \equiv 1728$ (mod $l$).  As above, there are at least two roots of $p_d(x) \equiv 0$ in $\mathbb{F}_l$ corresponding to $j=1728$, when $l \equiv 3$ (mod $4$) (since $a \not \equiv -1/a$ (mod $l$)).  On the other hand, there cannot be more than two roots each for the polynomials $p_l(x)$ and $p_{4l}(x)$, by the last displayed equations, because $p_d(x)$ has no more than $2\textsf{h}(-d)$ distinct roots (mod $l$), by Theorem 3.1.  Thus, both polynomials have exactly $2$ roots corresponding to $j=1728$, and this proves that $p_d(x) \equiv 0$ has $2\textsf{h}(-d)$ distinct roots in $\mathbb{F}_l$ in all cases.  $\square$ \medskip

This fact allows us to prove \bigskip

\noindent {\bf Theorem 4.2.} {\it If $l \equiv 4$ (mod $5$) and, with the previous notation, the minimal polynomial of $\eta=r(w/5)$ over $K=\mathbb{Q}(\sqrt{-d})$ is $m_d(x)$, where $d=l$ or $d=4l$, then the prime divisor $\mathfrak{l}=(\sqrt{-l})$ of $K$ does not divide $\textrm{disc}(m_d(x))$, so that the powers of $\eta$ form an $\mathfrak{l}$-integral basis of $F_1/K$.} \medskip

\noindent {\it Proof} We have $p_d(x)=m_d(x) \overline{m}_d(x)$ over $K$, where the bar denotes complex conjugation.  Clearly, $m_d(x) \equiv \overline{m}_d(x)$ (mod $\mathfrak{l}$), so that $p_d(x) \equiv m_d(x)^2$ (mod $\mathfrak{l}$).  Thus, $m_d(x)$ has no multiple roots (mod $\mathfrak{l}$), by Theorem 4.1, and this proves the theorem. $\square$
\bigskip

We now prove the following theorem, which will be useful in Section 7. \bigskip

\noindent {\bf Theorem 4.3.} {\it As in Theorem 1.1, let $h(-l)$ denote the class number of the field $K=\mathbb{Q}(\sqrt{-l})$.}  \begin{enumerate}[a)]
\item {\it If $l \equiv 4$ mod $5$, the polynomial $\hat H_{5,l}(x)$ has no irreducible factors of the form $q(x) = x^2+ax-1$ over $\mathbb{F}_l$.}
\item {\it If $l \equiv 1$ mod $5$, the number of irreducible quadratic factors $q(x)=x^2+ax-1$ dividing $\hat H_{5,l}(x)$ over $\mathbb{F}_l$ is:}
\begin{enumerate}[i)]
\item {\it $h(-l)$, if $l \equiv 1$ mod $4$;}
\item {\it $4h(-l)$, if $l \equiv 3$ mod $8$;}
\item {\it $2h(-l)$, if $l \equiv 7$ mod $8$.}
\end{enumerate}
\end{enumerate}

\noindent {\it Proof.} Assume $q(x)=x^2+\tilde a x-1$ is an irreducible quadratic dividing $\hat H_{5,l}(x)$.  Let $b$ be a root of this polynomial over $\mathbb{F}_l$.  Then the curve $E_5$ with the parameter $b$ is supersingular, and $z = b-\frac{1}{b} = -\tilde a \in \mathbb{F}_l$, which implies that the $j$-invariant
$$j =\frac{(b^4-228b^3+494b^2+228b+1)^3}{b(1-11b-b^2)^5}=-\frac{(z^2-228z+496)^3}{(z+11)^5} \in \mathbb{F}_l.$$
This is the $j$-invariant of a curve isogenous to $E_5(b)$, which is denoted $E_{5,5}$ in \cite{mor}, \cite{mor2}.  (Recall the alternate formula for $\hat H_{5,l}(b)$ in Section 2.1.)  Moreover, both roots of $q(x)$ lead to the same supersingular $j$-invariant, since the roots of $q(x)$ are $b$ and $-1/b$.  Hence, $q(x)$ must divide $G(x,j)$ in $\mathbb{F}_l$.  \medskip

First assume that $l \equiv 4$ (mod $5$).  As we proved in Propositions 3.2 and 3.3, the roots of $q(x)$ have the form $M(\rho)^5$, where $\rho \in \mathbb{F}_l$ is one of four (or two) roots of $p_d(x)$ (with $d=l$ or $4l$), and $M = S^k A^i$ are left coset representatives of the subgroup $H$ in $G_{60}$.  Since the only fifth root of $-1$ in $\mathbb{F}_l$ is $-1$ itself, it suffices to show that
$$S^k A^i(\rho) \cdot (S^k A^i)^l(\rho) \not \equiv -1 \ \ \textrm{in} \ \mathbb{F}_l.$$
Letting $\sigma = (\zeta \rightarrow \zeta^2)$, we have that $\zeta^l = \zeta^4 = \sigma^2(\zeta)$, so that $\sigma^2$ agrees with the automorphism $(\alpha \rightarrow \alpha^l)$ of $\mathbb{F}_{l^2}$ fixing $\mathbb{F}_l$.  Hence,
\begin{align*}
(S^k A^i)^l & = (S^k A^i)^{\sigma^2} = S^{-k} (A^{\sigma^2})^i = S^{-k} ((A^{-1} U)^\sigma)^i\\
&=S^{-k} (A^{-\sigma}U)^i = S^{-k} (U A U)^i = S^{-k}UA^i U.
\end{align*}
It follows that for $z \in \mathbb{F}_l$, and $A^i(z) =(az+b)/(cz+d)$,
\begin{align*}
S^k A^i(z) \cdot (S^k A^i)^l(z) &= \zeta^k \frac{az+b}{cz+d} \cdot \zeta^{-k} \left(-\frac{c-dz}{a-bz}\right)\\
&=- \frac{az+b}{cz+d} \cdot \frac{c-dz}{a-bz}=A^i(z) \cdot (UA^i U)(z).
\end{align*}
This expression can equal $-1$ if and only if $UA^iU(z) = UA^i(z)$, which holds if and only if $U(z)=z$ in $\mathbb{F}_l$.  But this condition is equivalent to $z^2+1=0$, and as we have already seen, this never holds for a root of $\hat H_{5,l}(X)$ in $\mathbb{F}_l$.  This proves part a).  \medskip

Now assume $l \equiv 1$ mod $5$.  With $\mathfrak{l}=(\sqrt{-l})$ in $K = \mathbb{Q}(\sqrt{-l})$, consider the factorization of $m_d(x)$ in $R_K/\mathfrak{l}$, where $m_d(x)$ is the minimal polynomial of $\eta$ over $K$, as in Theorem 4.2.  I claim that $m_d(x)$ factors (mod $\mathfrak{l}$) as a product of $\textsf{h}(-d)$ distinct quadratic factors of the form $q(x)=x^2+\tilde a x-1$.  The argument is a modification of the argument we gave in the proof of Proposition 3.2.  Recall that $F_1 = \Sigma_{\wp_5'}$ or $F_1 = \Sigma_{\wp_5'} \Omega_2$, as in Section 3.  In either case, $F_1=\Sigma(\eta) = \Sigma(b)$, where $b=\eta^5=r(w/5)^5$ and $\Sigma$ is as in the proof of Proposition 3.2.  For any singular $j$-invariant $\mathfrak{j}_i$, which is a root of $H_{-d}(x)$, there is a conjugate of $\eta$ over $K$, say $\eta_i=\eta^{\tau_i}$, which is a root of $G(x^5,\mathfrak{j}_i) = 0$.  It follows that the minimal polynomial $\mu_i(x)$ of $\eta_i$ over the field $\Sigma = K(\mathfrak{j}_i)$ divides both polynomials $G(x^5,\mathfrak{j}_i)$ and $m_d(x)$.  On the other hand, the nontrivial automorphism of $F_1/\Sigma$ is $\tilde \sigma: = (\eta \rightarrow -1/\eta)$, so that the roots of each factor $\mu_i(x)$ are invariant under $\tilde \sigma$.  However, $\eta^{\tau_i \tilde \sigma} = \eta^{\tilde \sigma \tau_i} = -1/\eta^{\tau_i}$ implies that $\mu_i(x) = (x-\eta_i)(x+1/\eta_i) = x^2+\tilde a_i x-1$, where $\tilde a_i = 1/\eta_i-\eta_i \in \Sigma$.  If $\mathfrak{p}$ is a fixed prime divisor of $\mathfrak{l}$ in $\Sigma$, then $\mathfrak{p}$ has degree $1$ since $\mathfrak{l}$ is principal, and $\sqrt{-l} \equiv 1$ (mod $2$) in case $\Sigma = \Omega_2$.  It follows that $m_d(x)$ has, for each distinct $j_i \equiv \mathfrak{j}_i$ (mod $\mathfrak{p}$), a factor $\mu_i(x) = x^2+\tilde a_i x-1$ (mod $\mathfrak{p}$) in $\Sigma/\mathfrak{p} \cong R_K/\mathfrak{l}$.  Using the same argument with $F_2=\Sigma(\xi)$ in place of $F_1$ yields a second factor $\tilde \mu_i(x)$ dividing $G(x^5,\mathfrak{j}_i) = 0$ and $\overline{m}_d(x)$ (the minimal polynomial of $\xi$ over $K$).  As in Proposition 3.2, $\mu_i(x) \not \equiv \tilde \mu_i(x)$ (mod $\mathfrak{p}$), if the $j$-invariant corresponding to $z_i=\eta_i^5-1/\eta_i^5$ is not $0$ or $1728$; and $\tilde \mu_i(x) \mid m_d(x)$ (mod $\mathfrak{p}$) since $m_d(x) \equiv \overline{m}_d(x)$ (mod $\mathfrak{l}$).  Furthermore, $l$ divides $d$, so the $j$-invariants $j_i \equiv \mathfrak{j}_i$ (mod $\mathfrak{p}$) are all supersingular, where $j_i$ is given by the rational function in $z_i$ in the first paragraph of the proof.  Using \cite[Theorem 6.2]{mor}, we see that $\hat H_{5,l}(x)$ cannot have any linear factors (mod $l$), and therefore the polynomials $\mu_i(x), \tilde \mu_i(x)$ remain irreducible (mod $\mathfrak{p}$). \medskip

If $\left(\frac{-3}{l}\right) = -1$, then by the above argument, $m_d(x)$ (for $d = 4l$, see the remarks in the Introduction after equation (1.1)) has at least one factor $\mu(x)$ corresponding to $j \equiv 0$.  Since $p_d(x) \equiv m_d(x)^2$ (mod $\mathfrak{p}$) is invariant under the mapping $T$, then $m_d(x)$ is also fixed by $T$ (up to a constant factor).  By the factorization of $g_0(x)$ in Section 2.1, $\mu(x)$ would have to be congruent to one of the quadratic factors of $x^4 - 3x^3 - x^2 + 3x + 1$ (see below), say
$$\mu(x) =  x^2 + \frac{-3 + \sqrt{5}}{2}x - 1,$$
in which case
$$(2x+1+\sqrt{5})^2 \mu(T(x)) = 4\sqrt{5}(x^2 + \frac{-3 - \sqrt{5}}{2}x - 1) = 4\sqrt{5} \tilde \mu(x).$$
This shows that if $m_d(x)$ is divisible (mod $\mathfrak{p}$) by one factor $\mu(x)$ corresponding to $j \equiv 0$, then it is also divisible by the conjugate factor $\tilde \mu(x)$.  Using (1.1) and (1.2) now shows that $m_d(x)$ has least
$$2\left(\frac{\textsf{h}(-d)-1}{2}\right)+1 = \textsf{h}(-d)$$
distinct quadratic factors of the form $q(x) = x^2+ \tilde a x-1$ (mod $\mathfrak{l}$), if $l \equiv 3$ (mod $4$); and at least $2(\frac{\textsf{h}(-d)}{2}) = h(-d)$ distinct quadratic factors of this form (mod $\mathfrak{l}$), if $l \equiv 1$ (mod $4$).  Hence, $m_d(x)$ is a product of $\textsf{h}(-d)$ distinct quadratics of the form $q(x)$, in all cases, as claimed.  \medskip

In the case under consideration, for each irreducible polynomial of the form $q(x) = x^2 +a x - 1$ dividing $\hat H_{5,l}(x)$, there is exactly one factor $\tilde q(x) = x^2 + \tilde a x -1$ of the same form which divides the polynomial $q(x^5)$.  For if $\alpha^5, \beta^5 = \alpha^{5l}$ (with $\alpha^5 \beta^5 = -1$) are the roots of $q(x)$, then the roots of quadratic factors of $q(x^5)$ have the form $\zeta \alpha, \zeta \beta = \zeta \alpha^l$, where $\zeta \in \mathbb{F}_l$ is a $5$-th root of unity.  Then $(\alpha \beta)^5 = -1$ implies that there is a unique $\zeta$ for which $(\zeta \alpha)(\zeta \beta) = \zeta^2 \alpha \beta = -1$.  Hence, the factors of the form $x^2 +\tilde a x -1$ of $m_d(x)$ modulo $\mathfrak{p}$ are in $1-1$ correspondence with factors of $\hat H_{5,l}(x)$ of the same form.  \medskip

Thus, $\hat H_{5,l}(x)$ has at least $h(-l) +  \textsf{h}(-4l)$ factors of the form $q(x)$ when $l \equiv 3$ (mod $4$), and $h(-l)$ such factors when $l \equiv 1$ (mod $4$).  This yields the formulas in part b).  It remains to show that these are all the factors of the form $q(x)$ which occur in $\hat H_{5,l}(x)$.  \medskip

To do this we use an argument similar to the proof of Proposition 3.3.  We show that the the expression
$S^i A^k(z) \cdot (S^i A^k(z))^l \neq -1$, if $S^iA^k \notin H$, $j \neq 0, 1728$, and $z$ is a root of $p_d(x)$ (mod $l$).  (The factors corresponding to $j = 0, 1728$ are included in the factors which we enumerated above.)  In this case the $5$-th roots of unity lie in $\mathbb{F}_l$, and $z^l \equiv -1/z =U(z)$ (mod $l$); so we must show
$$S^i A^k(z) \cdot S^i A^kU(z) \neq -1 \ \textrm{in} \ \mathbb{F}_l.$$
First note that we may assume $k \not \equiv 0$ (mod $3$), since 
$$S^i(z)\cdot S^i\left(\frac{-1}{z}\right) = -\zeta^{2i} \neq -1$$
unless $i=0$.  Now the claim for $k=1$ follows as in the proof of Proposition 3, on noting the following factorizations:
\begin{align*}
A(z) \cdot AU(z)+1 & = \frac{(1-\zeta^2)(z^2+(1-\sqrt{5})z-1)}{(z-1-\zeta^4)(z+1+\zeta+\zeta^3)},\\
SA(z) \cdot SAU(z)+1 & = \frac{(1-\zeta^4)(z-\zeta^2-\zeta^3)(z+1+\zeta^2+\zeta^3)}{(z-1-\zeta^4)(z+1+\zeta+\zeta^3)},\\
S^2A(z) \cdot S^2 AU(z)+1 & = \frac{(1-\zeta)(z^2-\frac{(3+\sqrt{5})}{2}z-1)}{(z-1-\zeta^4)(z+1+\zeta+\zeta^3)},\\
S^3A(z) \cdot S^3 AU(z)+1 & = \frac{(1-\zeta^3)(z^2-\frac{(3-\sqrt{5})}{2}z-1)}{(z-1-\zeta^4)(z+1+\zeta+\zeta^3)},\\
S^4A(z) \cdot S^4 AU(z)+1 & = \frac{(2+4\zeta+\zeta^2+3\zeta^3)z}{(z-1-\zeta^4)(z+1+\zeta+\zeta^3)}.
\end{align*}
This proves what we need, since $z \notin \mathbb{F}_l$ and
\begin{align*}
&(z^2+(1-\sqrt{5})z-1)(z^2+(1+\sqrt{5})z-1) = z^4+2z^3-6z^2-2z+1,\\
&(z^2-\frac{(3+\sqrt{5})}{2}z-1)(z^2-\frac{(3-\sqrt{5})}{2}z-1)=z^4-3z^3-z^2+3z+1;
\end{align*}
and these factors correspond to $j=1728$ and $j=0$.  For $k=2$ we use the fact that $A^2 = A^\sigma U$ from Section 2.2.  Then
$$S^{2i}A^2(z) \cdot S^{2i}A^2U(z) +1= S^{2i}A^\sigma U(z) \cdot S^{2i} A^\sigma(z) +1= \left(S^iA(z) \cdot S^i AU(z) + 1\right)^\sigma,$$
where $\sigma$ is the result of replacing $\zeta$ by $\zeta^2$ and leaving $z$ unchanged.  This is an automorphism (over the rational function field $\mathbb{Q}(z)$) and simply replaces $\sqrt{5} = \zeta-\zeta^2-\zeta^3+\zeta^4$ by $\sqrt{5}^\sigma = -\sqrt{5}$ in the relevant factors above.  Hence, the same argument applies. \medskip

This proves that, for $j \neq 0, 1728$, the only factors of $G(x^5,j)$ of the form $q(x)=x^2+\tilde a x-1$ over $\mathbb{F}_l$ are the factors of $p_d(x)$, for $d=l$ (when $l \equiv 3$ mod $4$) and $d=4l$.  Moreover, we know that $G(x,0) = g_0(x)^3$ and $G(x,1728) = g_{1728}(x)^2$, where $g_0$ and $g_{1728}$ are the polynomials in Section 2.1, and for these polynomials there are no factors (mod $l$) other than the ones guaranteed by Proposition 2.1c) and d).  This completes the proof.  $\square$. \bigskip

\noindent {\bf Corollary 4.4.} {\it If $l$ is a prime with $l \equiv 1$ (mod $5$), the polynomial $p_d(x)$ ($d=l$ or $4l$) is the product of the squares of $\textsf{h}(-d)$ distinct, irreducible, quadratic factors modulo $l$.  In particular, Theorem 4.2 also holds for primes $l \equiv 1$ (mod $5$).}

\noindent {\it Proof.} This is immediate from the above proof and $p_d(x) = m_d(x) \overline{m}_d(x) \equiv m_d(x)^2$ (mod $\mathfrak{l}$). 
$\square$

\section{The Hasse invariant for $E_7$.}

The Tate normal form for a point of order $7$ is the curve
$$E_7: \ \ Y^2+(1+d-d^2)XY+(d^2-d^3)Y=X^3+(d^2-d^3)X^2,$$
whose $j$-invariant is
$$j_7(d) = \frac{(d^2-d+1)^3(d^6-11d^5+30d^4-15d^3-10d^2+5d+1)^3}{d^7(d-1)^7(d^3-8d^2+5d+1)},$$
(see \cite{du}); and whose Hasse invariant is the value at $x=d$ of the polynomial
\begin{align*}
\hat H_{7,l}(x) & = (x^2-x+1)^r(x^6-11x^5+30x^4-15x^3-10x^2+5x+1)^r\\
& \times (x^{12}-18x^{11}+117x^{10}-354x^9+570x^8-486x^7+273x^6\\
& -222x^5+174x^4-46x^3-15x^2+6x+1)^s\\
& \times x^{7n_l}(x-1)^{7n_l}(x^3-8x^2+5x+1)^{n_l}J_l(j_7(x)).
\end{align*}
From \cite[Theorem 6.1]{mor} we have the following statements about the supersingular parameters for the curve $E_7$ over $\mathbb{F}_l$. \bigskip

\noindent {\bf Proposition 5.1.} {\it Assume $E_7$ is supersingular over $\mathbb{F}_l$, for a given prime $l \neq 2, 3, 7$.  Putting $a = 1+d-d^2$ and $b=d^2-d^3$, we have:}
\begin{enumerate}[i)]
\item {\it If $l \equiv 1$ (mod $7$), then $\mathbb{F}_l(a,b) = \mathbb{F}_{l^2}$.}
\item {\it If $l \equiv 6$ (mod $7$), then $\mathbb{F}_l(a,b) = \mathbb{F}_l$ or $\mathbb{F}_{l^2}$.}
\item {\it If $l \equiv 2,4$ (mod $7$), then $\mathbb{F}_l(a,b) = \mathbb{F}_{l^2}$ or $\mathbb{F}_{l^6}$.}
\item {\it If $l \equiv 3,5$ (mod $7$), then $\mathbb{F}_l(a,b) = \mathbb{F}_{l^2}, \mathbb{F}_{l^3}$ or $\mathbb{F}_{l^6}$.}
\end{enumerate}
\medskip

\noindent {\it Proof.} The proof in \cite[Theorem 6.1(ii), (iii)]{mor} is correct in case the $j$-invariant of $E_7$ satisfies $j \neq 0, 1728$.  However, if $j=0$ or $1728$, the proof given there is incomplete.  In the first paragraph of the proof marked ``Proof of (ii) and (iii)", the meromorphism $\pi(x,y) = (x^q, y^q)$ satisfies $\pi^2 = \zeta q$, where $\zeta$ is a root of unity in $\textsf{R}=End(E)$, for the supersingular elliptic curve
$$E: \ Y^2+aXY+bY=X^3+bX^2,$$
and $q=p^r$ is the smallest power of $p$ for which $a, b \in \mathbb{F}_q$.  The meromorphism $\pi$ must satisfy an equation of degree $2$ or less over $\mathbb{Z}$ in $\textsf{R}$.  From this it is concluded that $\zeta= \pm 1$.  (The same argument is given in \cite[p. 263]{h5}.)  However, it is possible for $\zeta$ to be a $4$-th or $6$-th root of unity and for $\pi$ still to be the root of a quadratic equation. \medskip

For example, suppose that $j(E_7)=0$ in $\mathbb{F}_p$.  Consider the possibilities $\pi^2 = \zeta q$ and $\zeta = \delta \omega$, where $\delta = \pm 1$ and $\omega$ is a primitive cube root of unity.  In this case we have
$$\pi^4 +\delta q \pi^2 + q^2 = (\omega^2+\omega+1)q^2 = 0;$$
and it is possible for the polynomial $x^4+\delta qx^2+q^2$ to factor.  If it does factor, then for some $c \in \mathbb{Z}$,
$$x^4+\delta qx^2+q^2=(x^2+cx\pm q)(x^2-cx\pm q)=(x^2 \pm q)^2-c^2x^2;$$
which implies that $\delta q = \pm 2q-c^2$, or $c^2 = \pm 2q-\delta q$ in $\mathbb{Z}$.  This implies that the plus sign holds in the first term.  If $\delta=1$, then $c^2=q=p^r$, implying that $r$ is even, and $\pi = \pm \omega^2 p^{r/2}$ in $\textsf{R}$.  In this case it is not hard to show as in \cite{mor} that $r \mid 2m$ and $2m \mid 3r$, where $m$ is the order of $p$ in the quotient group $\mathbb{Z}_n/\{\pm1\}$.  If $r=2m$ we obtain the contradiction $\pi=\pm p^m$ on using the proof of Theorem 6.1(i).  We obtain, finally, that $2m=3r$ in the case $\delta = +1$, so $r=2m/3$ and $\mathbb{F}_p(a, b) = \mathbb{F}_{p^{2m/3}}$.  If $\delta=-1$, then $c^2=3q$, implying that $q=3^{2s-1}$ is an odd power of $p=3$.  \medskip

Similarly, if $j(E_7)=1728$, then $\pi^2 = iq$ can lead to a quadratic equation for $\pi$.  This is because $\pi^4 = -q^2$, and $x^4+q^2$ factors if and only if $q^2 = 4c^4$ for some integer $c$ (by Capelli's theorem, \cite[p. 87]{co2}).  This implies that $q = 2c^2$, so $q=2^{2s-1}$ is an odd power of $2$.  Hence, the proof in \cite{mor} works in this case if we exclude the prime $p=2$.\medskip

These arguments show that when $p \neq 2,3$, $j(E)=0$, and $3 \mid m$, we have to add the possibility that $\mathbb{F}_p(a,b)=\mathbb{F}_{p^{2m/3}}$ to \cite{mor}, Theorem 6.1 (ii), (iii).  Thus, we must add the field $\mathbb{F}_{l^2}$ to the list of possibilities in parts (iii), (iv) of this theorem, since $m=3$ in these cases.
$\square$  \bigskip

The above considerations explain why it is possible to have a quadratic factor in parts (iii) and (iv) of the following theorem.\bigskip

\noindent {\bf Theorem 5.2.} {\it If $l \neq 2, 3, 7$ is prime, the irreducible factors of $\hat H_{7,l}(x)$ over $\mathbb{F}_l$ are:}
\begin{enumerate}[i)]
\item {\it quadratic, if $l \equiv 1$ (mod $7$);}
\item {\it linear or quadratic, if $l \equiv 6$ (mod $7$);}
\item {\it $x^2-x+1$ \ or \ sextic, if $l \equiv 2,4$ (mod $7$);}
\item {\it $x^2-x+1$, cubic or sextic, if $l \equiv 3,5$ (mod $7$).}
\end{enumerate}
\medskip

\noindent {\it Proof.} From the identity $d(a-1)=b$ and $a-1=d(1-d) \neq 0$ it is clear that $\mathbb{F}_l(a,b)=\mathbb{F}_l(d)$.  Thus the statements of (i)-(iv) follow directly from Proposition 5.1, for the factors of $\hat H_{7,l}(x)$ coming from the $12$-th degree factor (corresponding to $j(E_7)=1728$) and from $J_l(j_7(x))$.  On the other hand, it is clear that $x^2-x+1$ divides $\hat H_{7,l}(x)$ whenever $l \equiv 2$ (mod $3$).
We need to show that the factors over $\mathbb{F}_l$ of the polynomial
$$f_1(x) = x^6-11x^5+30x^4-15x^3-10x^2+5x+1,$$
conform to the statements in (i)-(iv), {\it when} when it divides the Hasse invariant.  \smallskip

A root of the polynomial $f_1(x)$ generates the real subfield $L_{21}^+ = \mathbb{Q}(\zeta_{21}+\zeta_{21}^{-1})$ of the $21$-st roots of unity, since one of its roots is
$$\alpha=\rho^5-\rho^4-5\rho^3+4\rho^2+4\rho, \ \ \rho = \zeta_{21}+\zeta_{21}^{20}.$$
Its discriminant is $\textrm{disc}(f_1(x)) = 2^{12} 3^3 7^5$.  Its Galois group satisfies $G(L_{21}^+/\mathbb{Q}) \cong \mathbb{Z}_{21}^\times/\{\pm 1\}$ and it is the class field over $\mathbb{Q}$ corresponding to the congruence group $\textsf{H}_1=\{1, 20\}$ (mod $21$).  Note that $f_1(x)$ is a factor of $\hat H_{7,l}(x)$ if and only if $l \equiv 2$ (mod $3$).  Hence we have
\begin{align*}
& l \equiv 1 \ (\textrm{mod} \ 7) \ \Rightarrow \ l \equiv 8 \ (\textrm{mod} \ 21), \ \textrm{with} \ 8^2 \equiv 1 \ (\textrm{mod} \  21);\\
& l \equiv 6 \ (\textrm{mod} \ 7) \ \Rightarrow \ l \equiv 20 \ (\textrm{mod} \ 21);\\
& l \equiv 2, 4 \ (\textrm{mod} \ 7) \ \Rightarrow \ l \equiv 2, 11 \ (\textrm{mod} \ 21), \ \textrm{with} \ 2^6 \equiv 11^6 \equiv 1 \ (\textrm{mod} \  21);\\
& l \equiv 3, 5 \ (\textrm{mod} \ 7) \ \Rightarrow \ l \equiv 17, 5 \ (\textrm{mod} \ 21), \ \textrm{with} \ 17^3 \equiv 5^3 \equiv 20 \ (\textrm{mod} \  21).
\end{align*}
The order of the residue class of $l$ in $\mathbb{Z}_{21}^\times/\{\pm 1\}$ determines the degree of the irreducible factors of $f_1(x)$ modulo $l$.  Therefore, the above congruences (mod $21$) show that (i)-(iv) do indeed hold for the degrees of these factors. $\square$  \bigskip

\noindent {\bf Theorem 5.3} {\it Let $\phi(x) = \frac{-1}{x-1}$. All irreducible cubic or sextic factors $f(x)$ of $\hat H_{7,l}(x)$ over $\mathbb{F}_l$ are invariant under the mapping $f(x) \rightarrow (x-1)^{\textrm{deg}(f)} f(\phi(x))$.  In addition:}\begin{enumerate}[a)]
\item{\it If $l \equiv 3, 5$ mod $7$, any root $d$ of $\hat H_{7,l}(x)$ satisfies $d^{l^2} = \phi(d)$.}
\item {\it All irreducible cubic factors of $\hat H_{7,l}(x)$ over $\mathbb{F}_l$ have the form
$$g_a(x)=x^3+ax^2-(a+3)x+1 \ (mod \ l).$$} 
\item{\it If $l \equiv 2, 4$ mod $7$, any root $d$ of $\hat H_{7,l}(x)$ satisfies $d^{l^2} = \phi^2(d)$.}
\item{\it All irreducible sextic factors $f(x)$ of $\hat H_{7,l}(x)$ over $\mathbb{F}_l$ have the form $f(x)=g_a(x) g_b(x)$, for conjugate elements $a, b \in \mathbb{F}_{l^2}$.}
\end{enumerate}
\medskip

\noindent {\it Proof.} a) We will make use of the meromorphism $\mu=l$ of the function field $\textsf{K}=\overline{\mathbb{F}}_l(x,y)$ of the supersingular curve $E_7$, as in \cite[Section 2]{brm} and \cite[Section 6]{mor}.  (Also see \cite{h6}, \cite{d} and \cite{ro}.)  Let $P_1=(0,0)$, and let $P_i=[i]P$ be the points of order $7$ on $E_7$ lying in $\langle P_1 \rangle$:
\begin{align*}
& P_1=(0,0), \ \ P_2=(d^3-d^2,d^3(d-1)^2), \ \ P_3 = (d^2-d, d(d-1)^2),\\
&P_4 = (d^2-d,d^2(d-1)^2), \ \ P_5 = (d^3-d^2,0), \ \ P_6 = (0, d^3-d^2).
\end{align*}
Further, let ${\bf p}_i$ represent the prime divisor of $\textsf{K}$ corresponding to the point $P_i$, and $\bf o$ the prime at infinity, so that
$$(x) = \frac{{\bf p}_1{\bf p}_6}{{\bf o}^2}, \ \ (y) = \frac{{\bf p}_1^2{\bf p}_5}{{\bf o}^3}.$$
As in \cite{mor}, we know that the endomorphism ${\bf q} \rightarrow \mu {\bf q} = l {\bf q}$ on prime divisors $\bf q$ of $\textsf{K}$ is given by $\mu {\bf q} = N_\mu(\bf q)^{\mu^{-1}}$, where $N_\mu$ is the norm from $\textsf{K}$ to $\textsf{K}^\mu$ and $\bf q$ is identified with the point $(x,y) \equiv (a, b)$ (mod $\bf q$) on $E_7$.  (Also see \cite{d}.)  It follows that
$${\bf p}_i^\mu=(\mu {\bf p}_{l^{-1}i})^\mu=N_\mu({\bf p}_{l^{-1}i}) = ({\bf p}_{l^{-1}i})^{l^2}, \ \ l^{-1} l \equiv 1 \ (\textrm{mod} \ 7),$$
since $\textsf{K}/\textsf{K}^\mu$ is purely inseparable of degree $[\textsf{K}:\textsf{K}^\mu]=l^2$.
If $l \equiv 3$ (mod $7$), $l^{-1} \equiv 5$ (mod $l$), so this implies that
$$(x^\mu) = \left(\frac{{\bf p}_5 {\bf p}_2}{{\bf o}^2}\right)^{l^2}, \ \ (y^\mu) = \left(\frac{{\bf p}_5^2 {\bf p}_4}{{\bf o}^3}\right)^{l^2}.$$
This gives in turn that
$$x^\mu = a(x-d^3+d^2)^{l^2}, \ \ y^\mu = c(dx+y-d^4+d^3)^{l^2}, \ \ a,c \in \overline{\mathbb{F}}_l, \eqno{(5.1)}$$
since $dx+y-d^4+d^3 = 0$ is the tangent to $E_7$ at $P_5$.  Applying $\mu$ to $dx+y-d^4+d^3$ gives
$$(dx+y-d^4+d^3)^\mu = \left(\frac{{\bf p}_4^2 {\bf p}_6}{{\bf o}^3}\right)^{l^2} = (d(d-2)x-y+d^3-d^2)^{l^2},$$
which implies that
$$dx^\mu+y^\mu-d^4+d^3 = e(d(d-2)x-y+d^3-d^2)^{l^2},$$
for some constant $e$.  Substituting from (5.1) and letting capital letters $X =x^{l^2}, Y=y^{l^2}$, etc., denote small letters raised to the power $l^2$, implies that
$$ad(X-D^3+D^2)+c(DX+Y-D^4+D^3)-d^4+d^3=e_1(D(D-2)X-Y+D^3-D^2).$$
The elements $1, X, Y$ are independent over $\overline{\mathbb{F}}_l$, so comparing coefficients of $Y$ on both sides gives $e_1=-c$.  Now comparing the coefficients of $X$ and the constant terms yields
\begin{align*}
ad+cD &=-cD(D-2)\\
ad(-D^3+D^2)+c(-D^4+D^3)-d^4+d^3 &=-c(D^3-D^2).
\end{align*}
Substituting for $ad+cD$ in the second equation from the first equation yields
$$-cD(D-2)(-D^3+D^2)-d^4+d^3=-c(D^3-D^2),$$
and therefore
$$c=\frac{d^4-d^3}{D^2(D-1)^3} , \ \ a =-\frac{d^3-d^2}{D(D-1)^2}.$$ \smallskip

On the other hand, applying the meromorphism $\mu$ to $(d(d-2)x-y+d^3-d^2)$ gives
\begin{align*}
(d(d-2)x-y+d^3-d^2)^\mu &= \left(\frac{{\bf p}_6^2 {\bf p}_2}{{\bf o}^3}\right)^{l^2},\\
d(d-2)x^\mu-y^\mu+d^3-d^2 &= e_2((d^2-d-1)x-y+d^3-d^2)^{l^2}, \ \textrm{or}
\end{align*}
\begin{align*}
ad(d-2)(X-D^3+D^2)&-c(DX+Y-D^4+D^3)+d^3-d^2\\
&=e_2((D^2-D-1)X-Y+D^3-D^2).
\end{align*}
Here we have $e_2=c$, so considering the coefficients of the $X$-terms gives
$$ad(d-2)-cD=c(D^2-D-1) \ \Longrightarrow \ \ -cD(D-1)(d-2)-cD=c(D^2-D-1);$$
this yields $(d-1)D^2-(d-2)D-1=(D-1)((d-1)D+1)=0$ and therefore
$$D=d^{l^2} = \frac{-1}{d-1}. \eqno{(5.2)}$$
In particular, $a=(d-1)^4$ and $c=-(d-1)^6$.  This shows that $-1/(d-1)$ is always a conjugate of $d$ in this case, and proves that the set of roots of any irreducible factor of $\hat H_{7,l}(x)$ is invariant under the map $x \rightarrow \phi(x) = -1/(x-1)$.  A similar argument works for the case $l \equiv 5$ (mod $7$). \medskip

\noindent b) The linear fractional map $\phi$ has order 3.  It follows that any cubic factor of $\hat H_{7,l}(x)$ must have the form
$$(x-a)(x-\phi(a))(x-\phi^2(a)) = x^3 +a_1x^2+a_2x+1,$$
where
$$a_1 = -\frac{a^3-3a+1}{a(a-1)} \ \ \textrm{and} \ \ a_2 = \frac{a^3-3a^2+1}{a(a-1)}.$$
Since $a_1+a_2 = \frac{3a-3a^2}{a(a-1)}= -3$, this proves the assertion of b).  Note that when $l \equiv 3, 5$ (mod $7$), the roots of an irreducible sixth degree factor $f(x)$ of $\hat H_{7,l}(x)$ also satisfy (5.2).  It follows that the cubic factors of $f(x)$ over $\mathbb{F}_{l^2}$ also have the form given in Part b). \medskip

\noindent c) When $l \equiv 2$ and $l^{-1} \equiv 4$ (mod $7$) we consider the effect of $\mu = l$ on $x, y$.  This time we have
$$(x^\mu) = \left(\frac{{\bf p}_4 {\bf p}_3}{{\bf o}^2}\right)^{l^2}, \ \ (y^\mu) = \left(\frac{{\bf p}_4^2 {\bf p}_6}{{\bf o}^3}\right)^{l^2},$$
which implies that
$$x^\mu = a(x+d-d^2)^{l^2}, \ \ y^\mu = c(d(d-2)x-y+d^3-d^2)^{l^2}, \ \ a,c \in \overline{\mathbb{F}}_l. \eqno{(5.3)}$$
Then
$$(d(d-2)x-y+d^3-d^2)^\mu = \left(\frac{{\bf p}_2^2 {\bf p}_3}{{\bf o}^3}\right)^{l^2}$$
yields the equation
$$d(d-2)x^\mu-y^\mu+d^3-d^2 = e_1((1-d^2)x+y+d^2(d-1)^2)^{l^2};$$
this gives
\begin{align*}
ad(d-2)(X+D-D^2)&-c(D(D-2)X-Y+D^3-D^2)+d^3-d^2\\
&=e_1((1-D^2)X+Y+D^2(D-1)^2).
\end{align*}
Thus, $e_1=c$ and the coefficients of $X$ satisfy
$$ad(d-2)-cD(D-2) =c(1-D^2) \ \Longrightarrow \ ad(d-2) = c(1-2D).\eqno{(5.4)}$$
Applying $\mu$ a third time to $(1-d^2)x+y+d^2(d-1)^2$ gives
$$((1-d^2)x+y+d^2(d-1)^2)^\mu = \left(\frac{{\bf p}_1^2 {\bf p}_5}{{\bf o}^3}\right)^{l^2};$$
hence
$$a(1-d^2)(X+D-D^2)+c(D(D-2)X-Y+D^3-D^2)+d^2(d-1)^2 = e_2Y.$$
Now this gives
$$a(1-d^2)+cD(D-2) = 0 \ \ \Longrightarrow \ a(1-d^2) = -cD(D-2).$$
Combining this with (5.4), we find that
$$(1-d^2)(1-2D) = -d(d-2)D(D-2),$$
or
$$0=(d^2-2d)D^2+(4d-2)D+1-d^2 = (dD+1-d)((d-2)D+1+d).$$
Hence,
$$D=d^{l^2} = \frac{d-1}{d} \ \ \textrm{or} \ \ D=d^{l^2} = -\frac{d+1}{d-2}.$$
The first case gives $d^{l^2} = \phi^2(d)$, with $\phi(x) = -1/(x-1)$, as above.  The linear fractional map $\rho(x) = -(x+1)/(x-2)$ has order $6$ and $\rho^2(x)=\phi(x)$, so the second case gives $d^{l^4}=\frac{-1}{d-1}=\phi(d)$.  If $d$ is a root of an irreducible factor of degree $6$, this gives that $d^{l^2}=\phi^2(d)$ in either case.  Note that $\rho^3(x) = \frac{x-2}{2x-1}$ must have $d$ as a fixed point, since $d^{l^6}=d$, and the fixed points of $\rho^3(x)$ are just the sixth roots of unity.  \medskip

Now factoring an irreducible sextic factor $f(x)$ over $\mathbb{F}_{l^2}$ and using the argument in part b) implies part d).  This completes the proof.
$\square$  \bigskip

\noindent {\bf Remark.} Note that the roots of the polynomials
\begin{align*}
f_0(x) & = x^2-x-1,\\
f_1(x) & = x^6-11x^5+30x^4-15x^3-10x^2+5x+1,\\
f_2(x) & = x^{12}-18x^{11}+117x^{10}-354x^9+570x^8-486x^7+273x^6\\
&  \ \ \ -222x^5+174x^4-46x^3-15x^2+6x+1,\\
f_3(x) & = x^3-8x^2+5x+1,
\end{align*}
are invariant under the map $\phi(x) = -1/(x-1)$, as is the $j$-invariant of $E_7$.  
\bigskip

\noindent {\bf Conjecture 1.} {\it Let $l \neq 2, 3, 7$ be a prime.}
\begin{enumerate}[a)]
\item {\it If $l \equiv 3,5$ (mod $7$), the number of distinct irreducible cubics which divide $\hat H_{7,l}(x)$ is}
\begin{align*}
h(-l), \ \ \textrm{if} \ \ l \equiv 1 \ (\textrm{mod} \ 4);\\
\left(3-\left(\frac{2}{l}\right)\right)h(-l), \ \ \textrm{if} \ \ l \equiv 3 \ (\textrm{mod} \ 4).
\end{align*}
{\it This is twice the number of supersingular $j$-invariants which lie in $\mathbb{F}_l$.}
\item {\it If $l \equiv 6$ (mod $7$), the number of distinct linear factors which divide $\hat H_{7,l}(x)$ is}
\begin{align*}
3h(-l), \ \ \textrm{if} \ \ l \equiv 1 \ (\textrm{mod} \ 4);\\
3\left(3-\left(\frac{2}{l}\right)\right)h(-l), \ \ \textrm{if} \ \ l \equiv 3 \ (\textrm{mod} \ 4).
\end{align*}
{\it This is six times the number of supersingular $j$-invariants which lie in $\mathbb{F}_l$.}
\end{enumerate}
\bigskip

Theorem 1.1 and Conjecture 1b) suggest the following regarding the Hasse invariant of the Tate normal form for a point of order $p$. \bigskip

\noindent {\bf Conjecture 2.} {\it If $p>3$ is prime and $l \equiv p-1$ (mod $p$), the number of linear factors of the Hasse invariant $\hat H_{p,l}(x)$ of the Tate normal form $E_p$ over $\mathbb{F}_l$ is $(p-1)$ times the number of supersingular $j$-invariants which lie in $\mathbb{F}_l$.} \bigskip

Theorem 1.1 shows that this conjecture is true for the prime $p=5$.  It is also true for $p=3$, if $E_3: Y^2+\alpha XY+Y=X^3$ is the Deuring normal form, by \cite{brm}, Theorems 1(d) and 5.

\section{Supersingular invariants for the Fricke groups.}

In this section, we will apply the results of Sections 3 and 5 to give a partial answer to a question and conjecture of Nakaya \cite{na}.  \medskip

As in \cite{na}, let $ss_p^{(N*)}(X)$ denote the supersingular polynomial for the Fricke group $\Gamma_0^*(N)$ introduced by Koike and Sakai \cite{sa1}, \cite{sa2}.  Here we focus on the cases $N=2, 3, 5, 7$.  As in \cite{na}, there is, for each such $N$, a Hauptmodul $j_N^*(\tau)$ for the group $\Gamma_0^*(N)$, given by:
\begin{align*}
j_2^*(\tau) &= \left(\frac{\eta(\tau)}{\eta(2\tau)}\right)^{24}+128+4096\left(\frac{\eta(2\tau)}{\eta(\tau)}\right)^{24},\\
j_3^*(\tau) &= \left(\frac{\eta(\tau)}{\eta(3\tau)}\right)^{12}+54+729\left(\frac{\eta(3\tau)}{\eta(\tau)}\right)^{12},\\
j_5^*(\tau) &= \left(\frac{\eta(\tau)}{\eta(5\tau)}\right)^6+22+125\left(\frac{\eta(5\tau)}{\eta(\tau)}\right)^6,\\
j_7^*(\tau) &= \left(\frac{\eta(\tau)}{\eta(7\tau)}\right)^4+13+49\left(\frac{\eta(7\tau)}{\eta(\tau)}\right)^4;
\end{align*}
where $\eta(\tau)$ is the Dedekind $\eta$-function.  The modular functions $j(\tau)$ and $j(N\tau)$ satisfy, in each case, a quadratic relation
$$R_N(j(\tau),j_N^*(\tau)) = R_N(j(N\tau),j_N^*(\tau)) = 0,$$
with coefficients in $\mathbb{Z}[j_N^*(\tau)]$, where:
\begin{align*}
R_2(X,Y) &= X^2-X(Y^2-207Y+3456)+(Y+144)^3,\\
R_3(X,Y) &= X^2-XY(Y^2-126Y+2944)+Y(Y+192)^3,
\end{align*}
\begin{align*}
R_5(X,Y) &= X^2-X(Y^5-80Y^4+1890Y^3-12600Y^2+7776Y+3456)\\
& \ \ +(Y^2+216Y+144)^3,\\
R_7(X,Y) &= X^2-XY(Y^2-21Y+8)(Y^4-42Y^3+454Y^2-1008Y-1280)\\
& \ \ +Y^2(Y^2+224Y+448)^3.
\end{align*}
See \cite{sa2}, \cite{na}.  Then, for each $N$, Sakai \cite{sa2} defined the polynomial $ss_p^{(N*)}(X)$ as follows:
$$ss_p^{(N*)}(X) = \prod_{j_N^*}{(X-j_N^*)} \in \mathbb{F}_p[X],$$
where $j_N^*$ runs over the distinct roots of $R_N(j,j_N^*) \equiv 0$ in $\overline{\mathbb{F}}_p$ for those values of $j$ which are supersingular in characteristic $p$, i.e. the roots of $ss_p(X)=0$. \medskip

Nakaya \cite[section 4]{na} asked if the supersingular values $j_N^*$ lie in $\mathbb{F}_{p^2}$, for any prime $p$, when $N$ is an element of the set
$$\mathfrak{S} = \{2, 3, 5, 7, 11, 13, 17, 19, 23, 29, 31, 41, 47, 59, 71\},$$
which is the set of primes $q$ for which all supersingular $j$-invariants lie in $\mathbb{F}_q$.  The set $\mathfrak{S}$ is also, of course, the set of prime factors of the order of the Monster group $\mathbb{M}$.  \medskip

We give an affirmative answer to this question for $N \in  \{2, 3, 5, 7\}$. \bigskip

\noindent {\bf Theorem 6.1.} {\it If $N \in \{2, 3, 5, 7\}$, all the supersingular values $j_N^*$ in characteristic $p$ lie in $\mathbb{F}_{p^2}$ for $p \neq N$.} \medskip

\noindent {\it Proof.} The proof for the case $N=2$ indicates the general idea, which is to find a parametrization for each of the genus zero curves $R_N(X,Y) = 0$.  To do this, we will substitute for $X$ the $j$-invariant of the Legendre normal form, $E_2: Y^2=X(X-1)(X-\lambda)$, namely
$$X=j_2(\lambda)= \frac{2^8(\lambda^2-\lambda+1)^3}{(\lambda^2-\lambda)^2},$$
and then factor the resulting expression in $Y$.  For the polynomial $R_2(X,Y)$ we find that $R_2(j_2(\lambda),Y)$ factors as a product of three linear polynomials in $Y$, whose roots are
$$\frac{16(2\lambda-1)^4}{\lambda(\lambda-1)}, \ \ \frac{-16(\lambda-2)^4}{\lambda^2(\lambda-1)}, \ \ \frac{16(\lambda+1)^4}{\lambda(\lambda-1)^2}.$$
Hence, a parametrization of $R_2(X,Y)$ is, for example,
$$X = \frac{2^8(\lambda^2-\lambda+1)^3}{(\lambda^2-\lambda)^2}, \ \ Y = \frac{16(\lambda+1)^4}{\lambda(\lambda-1)^2}, \ \ \ (N=2).$$
Any of the above three roots yields an equivalent parametrization.  If $j_2(\lambda)$ is supersingular, then from \cite[Prop. 1]{brm} we know that the corresponding value $\lambda$ lies in $\mathbb{F}_{p^2}$.  Hence, $Y=j_2^*$ lies in $\mathbb{F}_{p^2}$ when $j_2(\lambda)$ is supersingular in characteristic $p \neq 2$.  This proves the assertion for $N=2$.  \medskip

For $N=3$, we use the expression for the $j$-invariant of the Deuring normal form $E_3: Y^2+\alpha XY+Y=X^3$:
$$j_3(\alpha) = \frac{\alpha^3(\alpha^3-24)^3}{\alpha^3-27}.$$
Taking the root of the linear factor of $R_3(j_3(\alpha),Y)=0$ gives the parametrization
$$X = \frac{\alpha^3(\alpha^3-24)^3}{\alpha^3-27}, \ \ Y = \frac{\alpha^6}{\alpha^3-27}, \ \ \ (N=3).$$
Here we could also use the parametrization
$$X= \frac{\alpha^3(\alpha^3+216)^3}{(\alpha^3-27)^3}, \ \ Y = \frac{\alpha^6}{\alpha^3-27},$$
which arises from an isogenous curve to $E_3$ (\cite[eq. (2.5)]{mor5}); and this yields a number of equivalent parametrizations on substituting for $\alpha$ using the linear fractional maps in the group
$$G_{12} = \langle S, T \rangle, \ \ S(\alpha) = \omega \alpha, \ \ T(\alpha) = \frac{3(\alpha+6)}{\alpha-3}.$$
See \cite{mor5}.  If $j_3(\alpha)$ is supersingular in characteristic $p \neq 3$, then $\alpha \in \mathbb{F}_{p^2}$, by \cite[Thm. 2.3]{mor5}.  Hence $Y = j_3^* \in \mathbb{F}_{p^2}$, as well.  \medskip

For the case $N=5$, we substitute in $R_5(X,Y)$ for $X$ the $j$-invariant of the Tate normal form $E_5$, expressed in the following form:
$$j_5(b) = -\frac{(z^2+12z+16)^3}{z+11}, \ \ z = b-\frac{1}{b}.$$
This expression for $j_5(b)$ arises from the fact that it is invariant under the map $b \rightarrow -1/b$.  (See \cite{mor2}.)  Factoring $R_5(j_5(b),Y)$ yields the parametrization
$$X = -\frac{(z^2+12z+16)^3}{z+11}, \ \ Y = -\frac{z^2+4}{z+11}, \ \ \ (N=5).\eqno{(6.1)}$$
Now, from \cite[Thm. 6.2]{mor} we know that the factors of the Hasse invariant $\hat H_{5,p}(b)$ over $\mathbb{F}_p$ are linear, quadratic, or quartic in $b$.  If $b$ is the root of a linear or quadratic factor of $\hat H_{5,p}(b)$, then certainly $Y = j_5^* \in \mathbb{F}_{p^2}$.  On the other hand, the quartic factors of $\hat H_{5,p}(X)$ have the form $X^4+aX^3+cX^2-aX+1$, i.e. have roots that are invariant under $b \rightarrow -1/b$.  If $b$ is the root of such a factor, then $z=b-\frac{1}{b}$ satisfies the quadratic equation
$$z^2+az+c+2 = 0, \ \ a, c \in \mathbb{F}_p.$$
Therefore, if $j_5(b)$ is supersingular in characteristic $p \neq 5$, then $z \in \mathbb{F}_{p^2}$, from which it follows that $Y = j_5^* \in \mathbb{F}_{p^2}$, as well. \medskip

Now consider the case $N=7$.  Putting $X=j_7(d)$, the $j$-invariant of the Tate normal form $E_7$, into $R_7(X,Y)$ leads to the following parametrization:
\begin{align*}
X =& \frac{(d^2-d+1)^3(d^6-11d^5+30d^4-15d^3-10d^2+5d+1)^3}{d^7(d-1)^7(d^3-8d^2+5d+1)},\\
Y = & \frac{(d^2-d+1)^3}{d(d-1)(d^3-8d^2+5d+1)}, \ \ \ (N=7).
\end{align*}
Now we use Theorems 5.2 and 5.3.  First assume $p \neq 2,3,7$.  If $j_7(d)$ is supersingular in characteristic $p$, then $d$ is a root of $\hat H_{7,p}(X)=0$.  If $p \equiv \pm 1$ (mod $7$), then $d \in \mathbb{F}_{p^2}$ by Theorem 5.2i) and ii), which implies $j_7^*=Y \in \mathbb{F}_{p^2}$.  If $p \equiv 2, 3,4, 5$ (mod $7$), we know $d$ is the root of a cubic or sextic factor over $\mathbb{F}_p$ of a special form (excluding the trivial case when $d$ is a root of $x^2-x+1$).  Now form the resultant
\begin{align*}
\textrm{Res}_x((x^2-x+1)^3- &Yx(x-1)(x^3-8x^2+5x+1), x^3+ax^2-(a+3)x+1)\\
&= ((a+8)Y+a^2+3a+9)^3.
\end{align*}
This shows that when $d$ is the root of a cubic of the given form, over $\mathbb{F}_{p^2}$, then any corresponding value $j_7^*$ in characteristic $p \neq 7$ satisfies
$$j_7^* = Y = -\frac{a^2+3a+9}{a+8}.$$
This shows that $j_7^*$ lies in $\mathbb{F}_{p^2}$ in all cases, whether $a \in \mathbb{F}_p$ and $d$ is the root of a cubic factor, or $a \in \mathbb{F}_{p^2}$ and $d$ is the root of a sextic factor, which factors as a product of two cubics, by Theorem 5.3d).  \medskip

The statement for $N=7$ and $p=2,3$ follows by direct calculation, considering
$$R_7(0,Y) = Y^2(Y^2+224Y+448)^3,$$
since $j=0$ is the only supersingular invariant in characteristic $2$ or $3$.  This completes the proof of the theorem. $\square$

\section{Proof of Nakaya's conjecture for $N=5$.}

In the following theorem, we let $h(-5l)$ denote the class number of the quadratic field $\mathbb{Q}(\sqrt{-5l})$.  This theorem is analogous to theorems proved in \cite{mor4} and \cite{mor5} for quadratic factors of certain Legendre polynomials.  In those theorems the class numbers $h(-2p)$ and $h(-3p)$ show up in counting binomial quadratic factors of $P_{(p-e)/4}(x)$ and $P_{(p-e)/3}(x)$ (mod $p$).  In the following result, the class number $h(-5l)$ shows up in counting certain factors of $\hat H_{5,l}(x)$ (mod $l$).  These statements are also higher degree analogues of Theorems 1.1 and 4.3. \bigskip

\noindent {\bf Theorem 1.4.} \begin{enumerate}[a)]
\item {\it Assume $l \equiv 2,3$ mod $5$.  The number of irreducible quartics of the form $f(x) = x^4+ax^3+(11a+2)x^2-ax+1$ dividing $\hat H_{5,l}(x)$ over $\mathbb{F}_l$ is: }
\begin{enumerate}[i)]
\item {\it $\frac{1}{4}h(-5l)$, if $l \equiv 1$ mod $4$;}
\item {\it $\frac{1}{2}h(-5l)-1$, if $l \equiv 3$ mod $8$;}
\item {\it $h(-5l)-1$, if $l \equiv 7$ mod $8$.}
\end{enumerate}
\item {\it If $l \equiv 4$ mod $5$, the number of irreducible quadratic factors $k(x)=x^2+rx+s$ dividing $\hat H_{5,l}(x)$ over $\mathbb{F}_l$, with $r = \varepsilon ^5 (s-1)$ or $r = \bar \varepsilon^5 (s-1)$ in $\mathbb{F}_l$ is:}
\begin{enumerate}[i)]
\item {\it $\frac{1}{2}h(-5l)$, if $l \equiv 1$ mod $4$;}
\item {\it $h(-5l)-3$, if $l \equiv 3$ mod $8$;}
\item {\it $2h(-5l)-3$, if $l \equiv 7$ mod $8$.}
\end{enumerate}
\item {\it If $l \equiv 1$ mod $5$, the number of irreducible quadratic factors $k(x)=x^2+rx+s$ dividing $\hat H_{5,l}(x)$ over $\mathbb{F}_l$, with $r = \varepsilon ^5 (s-1)$ or $r = \bar \varepsilon^5 (s-1)$ in $\mathbb{F}_l$ is:}
\begin{enumerate}[i)]
\item {\it $\frac{1}{2}h(-5l)$, if $l \equiv 1$ mod $4$;}
\item {\it $h(-5l)-1$, if $l \equiv 3$ mod $8$;}
\item {\it $2h(-5l)-1$, if $l \equiv 7$ mod $8$.}
\end{enumerate}
\end{enumerate}
\bigskip

If the roots of $g(x) = x^2+ax+b=0$ are written as $\alpha^5, \beta^5$, then the relation $a = \varepsilon ^5 (b-1)$ becomes
$$\alpha^5+\beta^5 = \varepsilon^5 (1-\alpha^5 \beta^5).$$
Hence, the roots of a quadratic satisfying this condition yield a point $(X,Y) = (\alpha, \beta)$ on the curve $\mathcal{C}_5$:
$$\mathcal{C}_5: \ X^5 + Y^5  = \varepsilon^5 (1-X^5Y^5), \ \ \varepsilon = \frac{-1+\sqrt{5}}{2}.$$
Solutions to this equation in certain class fields of imaginary quadratic fields were exhibited in \cite{mor2}.
\bigskip

I will now show that Theorem 1.4 implies Nakaya's Conjecture 5 in \cite{na} for $N=5$.  In this conjecture $L(p)=S(\mathbb{F}_p)$ denotes the number of supersingular $j$-invariants in characteristic $p$ which lie in the prime field $\mathbb{F}_p$, as in \cite{na}.  Thus, $L(p)$ is the number of linear factors of the supersingular polynomial $ss_p(X)$ over $\mathbb{F}_p$. \bigskip

\noindent {\bf Theorem 7.1.} {\it Theorem 1.4 implies that the number of linear factors of the polynomial $ss_p^{(5*)}(X)$ over $\mathbb{F}_p$ is given by the formula}
\begin{align*}
L^{(5*)}(p) =  & \ \frac{1}{2}\left(1+\left(\frac{-p}{5}\right)\right)L(p)\\
& \ + \frac{1}{8}\Big\{2+\left(1-\left(\frac{-1}{5p}\right)\right)\left(2+\left(\frac{-2}{5p}\right)\right)\Big\}h(-5p)\\
= & \ \frac{1}{4}\left(1+\left(\frac{p}{5}\right)\right)h(-p)+\frac{1}{4} h(-5p), \ \textrm{if} \ p \equiv 1 \ \textrm{mod} \ 4;\\
= & \ \left(1+\left(\frac{p}{5}\right)\right)h(-p)+\frac{1}{2} h(-5p), \ \textrm{if} \ p \equiv 3 \ \textrm{mod} \ 8;\\
= & \ \frac{1}{2}\left(1+\left(\frac{p}{5}\right)\right)h(-p)+ h(-5p), \ \textrm{if} \ p \equiv 7 \ \textrm{mod} \ 8.
\end{align*}
\medskip

\noindent {\it Proof.} {\it Case 1: $p \equiv 2, 3$ mod $5$.} From \cite[Theorem 6.2]{mor} we know that the factors of $\hat H_{5,p}(x)$ are either $x^2+1$ (only when $p \equiv 3$ (mod $4$), since $x^2+1$ corresponds to the $j$-invariant $j=1728$) or irreducible quartic polynomials satisfying $x^4 f(-1/x) \equiv f(x)$ (mod $p$).  On the other hand, the proof of Theorem 6.1 shows that linear factors of $ss_p^{(5*)}(x)$ have roots of the form $a = -\frac{z^2+4}{z+11} \in \mathbb{F}_p$, where $z=b-\frac{1}{b}$ and $b$ is the parameter defining the curve $E_5$.  Making this substitution for $z$ shows that $b$ is a root of the polynomial
$$g_a(x) = x^4 + ax^3 + (11a+2)x^2 -ax+1.$$
Additionally, when the factor $x^2+1$ is present, $(x^2+1)^2 = x^4+2x^2+1$ has this form with $a=0$.  Theorem1.4a) implies that there are exactly $\frac{1}{4} h(-5p), \frac{1}{2} h(-5p)-1$ or $h(-5p)-1$ irreducible factors of this form which divide $\hat H_{5,p}(x)$, depending on the residue class of $p$ (mod $8$), giving this many values of $a \in \mathbb{F}_p$ which are roots of $ss_p^{(5*)}(x)$; and the factor $x^2+1$ in the case $p \equiv 3$ (mod $4$) yields the additional
value $a=0$.  This proves the formula for $L^{(5*)}(p)$ in this case. \medskip

\noindent {\it Case 2: $p \equiv 4$ mod $5$.} In this case, \cite[Theorem 6.2]{mor} says that $\hat H_{5,p}(x)$ has only linear and quadratic factors.  We count the number of distinct factors of the form $g(x)$ above which share roots with $\hat H_{5,p}(x)$.  Note that
$$d_a=\textrm{disc}(g_a(x)) =5^3a^2(a^2-44a-16)^2.$$
By the Pellet-Stickelberger-Voronoi Theorem (\cite[Appendix]{brm}), the number $n_p$ of irreducible factors of $g(x)$ mod $p$ satisfies
$$\left(\frac{d_a}{p}\right) = (-1)^{4+n_p}, \ \ \textrm{if} \ d_a \neq 0.$$
It follows that if $d_a \neq 0$, $g(x)$ has either $2$ or $4$ irreducible factors.  The polynomial $g(x)$ satisfies the formula $x^4 g(-1/x) = g(x)$, so that if $g(x)$ has one linear factor, it must have two, implying that it must split completely.  (A root $r$ in $\mathbb{F}_p$ cannot satisfy $r = -1/r$, since this would imply $r^2+1=0$ and $p \equiv 1$ mod $4$; but factors of $x^2+1$ cannot divide $\hat H_{5,p}(x)$ for such primes, since $p$ is not supersingular for $j(E_5(i)) = 1728$.  Alternatively, $g_a(x)$ only has roots in common with $x^2+1$ when $a=0$ and $g_0(x)=(x^2+1)^2$.)  Hence, $g_a(x)$ either splits completely or is a product of two quadratics. \medskip

We first count the linear factors of $ss_p^{(5*)}(x)$ coming from the linear factors of $\hat H_{5,p}(x)$.  By Theorem 1.1, the latter polynomial has $4L(p)$ linear factors, and by Theorem 1.2, four of these linear factors correspond to the same supersingular $j$-invariant.  The mapping taking
$$b \rightarrow a = -\frac{z^2+4}{z+11} = -\frac{(b^2+1)^2}{b(b^2+11b-1)} \eqno{(7.1)}$$
is a $4-1$ mapping, except for the values of $a$ for which $d_a=0$.  Now $a=0$ only corresponds to $b^2+1=0$, which never has roots in $\mathbb{F}_p$.  Since $\left(\frac{5}{p}\right)=+1$ in this case, the polynomial $a^2-44a-16$, whose discriminant is $2^4 5^3$, factors.  However, its roots must correspond to supersingular $j$-invariants, where
$$j=-\frac{(z^2+12z+16)^3}{z+11} \ \ \textrm{and} \ \ a^2-44a-16=\frac{(z^2+22z-4)^2}{(z+11)^2}=0.$$
Taking the resultant
\begin{align*}
\textrm{Res}_z(&(z^2+12z+16)^3+j(z+11),z^2+22z-4)\\
&= -5^3(j^2-1264000j-681472000) = -5^3 H_{-20}(j)
\end{align*}
shows that the corresponding $j$-invariants are roots of the polynomial $H_{-20}(X)$, and are therefore supersingular if and only if $\left(\frac{-20}{p}\right) = -1$, i.e., if and only if $p \equiv 3$ (mod $4$).  If $p \equiv 1$ mod $4$, it follows that the map (7.1) is $4-1$ everywhere, so there are exactly $L(p)$ linear factors of $ss_p^{(5*)}(x)$ arising from linear factors of $\hat H_{5,p}(x)$.  On the other hand, if $p \equiv 3$ (mod $4$), the quadratic $a^2-44a-16$ has two roots in $\mathbb{F}_p$, and the values of $b$ corresponding to these roots satisfy
$$b^4+22b^3-6b^2-22b+1 \equiv 0 \ (\textrm{mod} \ p).$$
The discriminant of this quartic is $2^8 5^9$, so it also has $4$ roots in $\mathbb{F}_p$ whenever it has one.  Furthermore, the roots of $x^4+22x^3-6x^2-22x+1$ in characteristic zero are real and lie in the real subfield $\Sigma^+$ of $\Sigma=\mathbb{Q}(\zeta_{20})$.  In this case, $p \equiv 19$ (mod $20$), so $p$ splits in $\Sigma^+$ and therefore the above congruence in $b$ has four distinct roots in $\mathbb{F}_p$.  This implies that for two values of $a \in \mathbb{F}_p$, $g_a(x) = (x-b_1)^2(x-b_2)^2$, where each $b_i \in \mathbb{F}_p$ is a root of $\hat H_{5,p}(x) = 0$ and $b_2=-1/b_1$.  Hence, there are altogether $\frac{4L(p)-4}{4}+2=L(p)+1$ linear factors of $ss_p^{(5*)}(x)$ coming from linear factors of the Hasse invariant, when $p \equiv 3$ (mod $4$). \medskip

Now consider values of $a \in \mathbb{F}_p$ coming from quadratic factors of $\hat H_{5,p}(x)$.  If $g_a(x)=k_1(x)k_2(x)$, then the quadratics $k_i(x)$ are either invariant under $\phi(x) = -1/x$ or are mapped into each other by $\phi$, meaning that $x^2k_1(\phi(x)) = ck_2(x)$, for some constant $c$.  If $k_1(x) = x^2+rx+s$ is invariant, then
$$x^2 k_1(-1/x) = 1-rx+sx^2 = s(x^2+rx+s),$$
implies that $r=0$ and $s=1$; or $-r = rs$, so $s= -1$.  Thus, either $k_1(x) = k_2(x) = x^2+1$ or $k_1(x)=x^2+rz-1$, with a similar expression for $k_2(x)$.  However, Theorem 4.3a) implies that none of the factors of $\hat H_{5,p}(x)$ have this form, for $p \equiv 4$ (mod $5$).  Hence, the $k_i(x)$ must be mapped to each other by $\phi$.  Now we have
\begin{align*}
g_a(x) &= \frac{1}{s}(x^2+rx+s)(1-rx+sx^2)\\
&=x^4+\left(r-\frac{r}{s}\right)x^3+\frac{s^2-r^2+1}{s}x^2-\left(r-\frac{r}{s}\right)x+1,
\end{align*}
which gives
$$a = r-\frac{r}{s} \ \ \textrm{and} \ \ 11a+2 = \frac{s^2-r^2+1}{s}.$$
Solving for $r$ is these two equations gives
$$r = \frac{-11+5\sqrt{5}}{2}(s-1) = \varepsilon^5 (s-1) \ \ \textrm{or} \ \ r = \frac{-11-5\sqrt{5}}{2}(s-1) = \bar \varepsilon^5 (s-1).$$
These are exactly the conditions in Theorem 1.4b).  If $p \equiv 1$ mod $4$, there are $\frac{1}{2}h(-5p)$ quadratic factors satisfying one of these conditions.  This yields $\frac{1}{4}h(-5p)$ pairs of quadratics satisfying the above conditions, and therefore this many linear factors of $ss_p^{(5*)}(x)$.  If $p \equiv 3$ (mod $4$), then $x^2+1$ satisfies these condtions with $r = 0$ and $s = 1$.  Subtracting $1$ from the number in Theorem 1.4b) yields either $\frac{1}{2}h(-5p)-2$ or $h(-5p)-2$ pairs of quadratics yielding values of $a$ in $\mathbb{F}_p$, and adding $a=0$ gives
$\frac{1}{2}h(-5p)-1$ or $h(-5p)-1$ total linear factors of $ss_p^{(5*)}(x)$ arising from quadratic factors of $\hat H_{5,p}(x)$.  \medskip

Now, \cite[Theorem 6.2]{mor} says that $\hat H_{5,p}(x)$ only has linear or quadratic factors when $p \equiv 4$ (mod $5$).  Hence all linear factors of $ss_p^{(5*)}(x)$ arise in one of the ways described above, and this yields
\begin{align*}
&L(p)+\frac{1}{4}h(-5p) \ \textrm{linear factors}, \ \textrm{if} \ p \equiv 1 \ (\textrm{mod} \ 4);\\
&(L(p)+1) + (\frac{1}{2}h(-5p)-1) = L(p) + \frac{1}{2}h(-5p) \ \textrm{linear factors}, \ \textrm{if} \ p \equiv 3 \ (\textrm{mod} \ 8);\\
&(L(p)+1) + (h(-5p)-1) = L(p) + h(-5p) \ \textrm{linear factors}, \ \textrm{if} \ p \equiv 7 \ (\textrm{mod} \ 8).
\end{align*}
This completes the proof of the theorem in this case. \medskip

\noindent {\it Case 3: $p \equiv 1$ mod $5$.} In this case, $\hat H_{5,p}(x)$ is only divisible by quadratic factors, by the results of \cite{mor}.  The analysis in Case 2 applies to this case, except that the polynomial $x^4+22x^3-6x^2-22x+1$ now factors into a product of two quadratics (mod $p$), since $p \equiv 3$ (mod $4$) and $p \equiv 1$ (mod $5$) imply that $p \equiv 11$ (mod $20$); so that $p$ has order $2$ in $\mathbb{Z}_{20}^\times/\langle \pm1 \rangle$.  This gives two quadratic polynomials $q(x) = x^2+rx-1$, which are factors of $\hat H_{5,p}(x)$, for which $q(x)^2 = g_a(x)$, with $a=2r$.  Together with $x^2+1$, there are three quadratics which yield values of $a \in \mathbb{F}_p$.  Moreover, factoring $x^4+22x^3-6x^2-22x+1$ over $\mathbb{Q}(\sqrt{5})$ gives
\begin{align*}
 & x^4+22x^3-6x^2-22x+1\\
 & = (x^2+(11+5\sqrt{5})x-1)(x^2+(11-5\sqrt{5})x-1)=q_1(x)q_2(x);
 \end{align*}
and both of these factors satisfy the condition that $r = \varepsilon^5 (s-1)$ or $r = \bar \varepsilon^5 (s-1)$.  Hence, these three quadratics are included in the counts of Theorem 1.4c).  The remaining quadratics yield either $\frac{1}{2}h(-5p)-2$ or $h(-5p)-2$ values of $a \in \mathbb{F}_p$ when $p \equiv 3$ (mod $4$).  Hence, there are altogether $\frac{1}{2}h(-5p)+1$ or $h(-5p)+1$ values of $a$ arising from quadratics satisfying the $r/s$ conditions. \medskip

Now we must account for the $a$'s which arise from invariant quadratics under $\phi$.  Since the above two polynomials $q_1(x)$ and $q_2(x)$ have already been counted, and have the form $x^2+rx-1$, Theorem 4.3b) shows that this leaves $2L(p)-2$ polynomials of this form.  I claim that these polynomials pair up in a unique way to give $L(p)-1$ values of $a \in \mathbb{F}_p$, when $p \equiv 3$ (mod $4$).  Let $k_1(x) = x^2+rx-1, k_2(x)=x^2+r'x-1$ be two polynomials of this form.  Then $k_1(x) k_2(x) = g_a(x)$ over $\mathbb{F}_p$ if and only if
$$x^4+(r+r')x^3+(rr'-2)x^2-(r+r')x+1=x^4+ax^3+(11a+2)x^2-ax+1,$$
which holds if and only if
\begin{align*}
r+r' & = a,\\
rr' & = 11a+4.
\end{align*}
It follows that $r$ is a solution of $r^2-ar+11a+4=0$, which gives that $a=\frac{r^2+4}{r-11}$, where $r \not \equiv 11$ (mod $p$).  This shows that there is at most one value of $a$ corresponding to a given factor $k_1(x)$.  Then the factor $k_2(x)$ is determined by $r'=\frac{11r+4}{r-11}$, since
$$x^2-ax+11a+4=\left(x-r\right)\left(x-\frac{11r+4}{r-11}\right),$$
so there is at most one factor $k_2(x)$ satisfying the requirement that $k_1(x) k_2(x) \equiv g_a(x)$ (mod $p$).  On the other hand, given the factor $k_1(x)=x^2+rx-1$ of $\hat H_{5,p}(x)$, we know $r \not \equiv 11$, since $x^2+11x-1$ is a factor of the discriminant $\Delta(x)$ (see Section 2).  By (6.1) with $z=-r$, there is a $j$-invariant
$$j(r)=\frac{(r^2-12r+16)^3}{r-11}$$
in $\mathbb{F}_p$ for which $j_5^*(r)=\frac{r^2+4}{r-11}=a$ is a corresponding root of $R_5(X,Y) \equiv 0$ (mod $p$).  Furthermore, $r'=\frac{11r+4}{r-11}$ satisfies $j_5^*(r')=a$ and
$$j(r') = j\left(\frac{11r+4}{r-11}\right) = \frac{(r^2+228r+496)^3}{(r-11)^5}.$$
From \cite[eq. (2.2b)]{mor2} $j(r')$ is the $j$-invariant of a curve $E_{5,5}$ which is $5$-isogenous to $E_5(b)$.  We just have to show that $j(r)$ is supersingular.  This is easy, because a root $b$ of $k_1(x)=0$ is a root of $\hat H_{5,p}(x) =0$ over $\mathbb{F}_p$, and the form of $k_1(x)$ yields $-r=b-\frac{1}{b}$, so that $j(r)=j(E_5(b))$ is certainly supersingular.  It follows that $j(r')$ is supersingular, as well.  Therefore, setting $-r'=b'-\frac{1}{b'}$ gives that $b'$ is a common root of $\hat H_{5,p}(x) =0$ and $k_2(x)=x^2+r'x-1$; the latter is an irreducible factor of $\hat H_{5,p}(x)$ over $\mathbb{F}_p$, since $r' \in \mathbb{F}_p$ and $\hat H_{5,p}(x)$ has no linear factors.  This shows that when $p \equiv 3$ (mod $4$), the $2L(p)-2$ factors of the form $k_1(x)$ pair up in a unique way to give $L(p)-1$ values of $a \in \mathbb{F}_p$.  \medskip

When $p \equiv 1$ (mod $4$), Theorem 4.3b) yields $2L(p)/2=L(p)$ pairs of invariant quadratics whose product yields a polynomial of the form $g_a(x)$, by the same arguments.  Hence, Theorem 1.4c) and Theorem 4.3b) yield
\begin{align*}
&L(p)+\frac{1}{4}h(-5p) \ \textrm{linear factors}, \ \textrm{if} \ p \equiv 1 \ (\textrm{mod} \ 4);\\
&(L(p)-1) + (\frac{1}{2}h(-5p)+1) = L(p) + \frac{1}{2}h(-5p) \ \textrm{linear factors}, \ \textrm{if} \ p \equiv 3 \ (\textrm{mod} \ 8);\\
&(L(p)-1) + (h(-5p)+1) = L(p) + h(-5p) \ \textrm{linear factors}, \ \textrm{if} \ p \equiv 7 \ (\textrm{mod} \ 8).
\end{align*}
This completes the proof.  $\square$

\medskip

\noindent Dept. of Mathematical Sciences, LD 270

\noindent Indiana University - Purdue University at Indianapolis (IUPUI)

\noindent Indianapolis, IN 46202

\noindent {\it e-mail: pmorton@iupui.edu}

\end{document}